\documentclass[12pt]{amsart}
\usepackage{amssymb,latexsym}  
\usepackage[all]{xy}

\newcommand{\charic}{{\operatorname{char}\,}} 
\newcommand{\GL}{\operatorname{GL}}
\newcommand{\adj}{\operatorname{adj}}
\newcommand{\sing}{\operatorname{sing}}
\newcommand{\isom}{{\, \cong \,}}
\newcommand{\Q}{{\mathbb Q}}
\newcommand{\F}{{\mathbb F}}
\newcommand{\K}{{\mathbb K}}
\newcommand{\Aff}{{\mathbb A}}
\newcommand{\PP}{{\mathbb P}}
\newcommand{\Gm}{{\mathbb G}_m}
\newcommand{\Z}{{\mathbb Z}}
\newcommand{\pr}{{\operatorname{pr}}}
\newcommand{\x}{{\bf x}}
\newcommand{\Kbar}{{\overline{K}}}
\newcommand{\Gal}{{\operatorname{Gal}}}
\newcommand{\ord}{{\operatorname{ord}}}
\newcommand{\SL}{{\operatorname{SL}}}
\newcommand{\PGL}{{\operatorname{PGL}}}
\newcommand{\Aut}{\operatorname{Aut}}
\newcommand{\Gr}{\operatorname{Gr}}
\newcommand{\Qbar}{{\overline{\Q}}}
\newcommand{\CL}{{\mathcal{L}}}
\newcommand{\G}{{\mathcal{G}}}
\newcommand{\FF}{{\mathcal{F}}}
\newcommand{\pf}{\operatorname{pf}}
\newcommand{\ra}{{\, \longrightarrow \,}}
\newcommand{\rank}{\operatorname{rank}}
\newcommand{\sign}{{\operatorname{sign}}} 
\newcommand{\la}{{\lambda}}
\newcommand{\ratto}{{ - \, \to}}

\newenvironment{Proof}{\par\noindent{\sc Proof:}}%
                      {\hspace*{\fill}\nobreak$\Box$\par\medskip}
\newenvironment{ProofOf}[1]{\par\noindent{\sc Proof of #1:}}%
                      {\hspace*{\fill}\nobreak$\Box$\par\medskip}

\newtheorem{Proposition}{Proposition}[section]
\newtheorem{Theorem}[Proposition]{Theorem}
\newtheorem{Lemma}[Proposition]{Lemma}
\newtheorem{Algorithm}[Proposition]{Algorithm}

\theoremstyle{definition}

\newtheorem{Definition}[Proposition]{Definition}
\newtheorem{Remark}[Proposition]{Remark}

\begin{document}

\date{9th October 2006}
\title[The invariants of a genus one curve]
{The invariants of a genus one curve} 
\author{Tom Fisher}
\address{University of Cambridge,
        DPMMS, Centre for Mathematical Sciences,
         Wilberforce Road, Cambridge CB3 0WB, UK}
 \email{T.A.Fisher@dpmms.cam.ac.uk}

\begin{abstract}
It was first pointed out by Weil \cite{Weil1954}
that we can use classical
invariant theory to compute the Jacobian of a genus one curve.
The invariants required for curves of degree $n=2,3,4$ were
already known to the nineteenth centuary invariant theorists. 
We have succeeded in extending
these methods to curves of degree $n=5$, where 
although the invariants are too large to write down as 
explicit polynomials, we have found a practical 
algorithm for evaluating them.
\end{abstract}

\maketitle

\section{Introduction}

We work throughout over a perfect field $K$ with algebraic closure~$\Kbar$.
In this introduction we further assume that $\charic(K) \not=2,3$.
Let $C$ be a smooth curve of genus one defined over $K$ and let 
$D$ be a $K$-rational divisor on $C$ of degree~$n$. 
If $n=1$ then $C(K) \not= \emptyset$, and 
$C$ is defined by a Weierstrass equation.
If $n \ge 2$ then the complete linear
system $|D|$ defines a morphism $C \to \PP^{n-1}$. 
If $n \ge 3$ then this morphism is an embedding and we call the
image a genus one normal curve of degree~$n$. 
For $n \le 5$ the pair $(C,D)$ is described by data of the following form.
\begin{Definition}
\label{def:1.1}
A genus one model of degree $n=1,2,3,4,5$ is \\
(i) if $n=1$ a Weierstrass equation \\
(ii) if $n=2$ a binary quartic \\
(iii) if $n=3$ a ternary cubic \\
(iv) if $n=4$ a pair of quadrics in $4$ variables \\
(v) if $n=5$ a $5 \times 5$ alternating matrix of 
linear forms in $5$ variables.
\end{Definition}
The equations defined by a genus one model of degree 5
are the $4 \times 4$ Pfaffians of the matrix.
It is a classical fact that every genus one normal quintic 
is defined by equations of this form.
An algorithm for computing these matrices is given in 
\cite{g1pf}, based on the Buchsbaum-Eisenbud structure
theorem \cite{BE1}, \cite{BE2} for Gorenstein ideals of codimension~3. 

We write $X_n$ for the (affine) space of all genus one models of degree~$n$
and $K[X_n]$ for its co-ordinate ring.
In \S\ref{invars} we specify a linear algebraic group $G_n$ acting on $X_n$.
In each case we find that the ring of invariants $K[X_n]^{G_n}$
is a polynomial ring in two variables. 
(This could be deduced by a theorem of Kempf \cite[Theorem 2.4]{Kempf}
if $\charic(K) = 0$.)
We label the generators $c_4$ and $c_6$.
In the case $n=1$ these are the usual polynomials defined in
\cite[Chapter III]{Si1}. In the cases
$n=2,3,4,5$ we find that $c_4$ and $c_6$ are homogeneous polynomials 
of degree $4n/(6-n)$ and $6n/(6-n)$. Moreover, as first pointed
out by Weil \cite{Weil1954} in the case $n=2$, the invariants $c_4$ and 
$c_6$ give a formula for the Jacobian. The invariants for $n=2,3,4$ have been 
known since the nineteenth century, and are surveyed in \cite{Mc+}.

Our work has two main goals. The first is to 
describe the ring of invariants % , and its connection with the Jacobian, 
in a way that emphasises the similarities between the cases $n=1,2,3,4,5$. 
The second is to give practical methods for evaluating the invariants.
In both instances our main original contribution is in the case $n=5$.

A relatively easy argument (reduction to the case $n=1$) shows that 
if invariants $c_4$ and $c_6$ of the expected degrees exist
then they generate the ring of invariants.
In the cases $n=2,3$ the existence of these invariants 
is settled by writing them down as explicit polynomials.
In the case $n=4$ there is a classical trick for reducing
to the case $n=2$. But in the case $n=5$ the invariants
are too large to write down as explicit polynomials.
This makes it difficult to show they exist.

One of the tools of classical invariant theory is the
so-called symbolic notation, as described in \cite{GY}.
This is an extremely compact notation for specifying invariants. 
For example in the case $n=3$ the invariants may be written
(see \cite[\S\S220,221]{Salmon} or \cite[\S4.5]{St})
\[ 
\begin{array}{rcl}
c_4 & = & 54 \times (abc)(bcd)(cda)(dab) \\
c_6 & = & 972 \times (abc)(abd)(bce)(caf)(def)^2.
\end{array}   \]
By introducing non-commuting symbols it is possible to write
down similar expressions in the case $n=5$.
But we have no way of showing these invariants are non-zero 
without expanding them as explicit polynomials.
As remarked above, this is not feasible.

In principle we could use the representation theory of Lie algebras,
specifically the Weyl character formula, to compute the dimension 
of the vector space of invariants of any given degree.
For details in the case $n=3$ we refer to
\cite[Exercise 13.20]{FH} or \cite[\S4.4]{St}. 
Unfortunately, when we tried this approach in the case $n=5$, 
we were again defeated by combinatorial explosion.

The plan of the paper is as follows. 
In \S\ref{geominvar} we explain the role played by the 
invariant differential in computing the 
Jacobian of a genus one curve.
In \S\ref{sec:g1models} we revisit and motivate our definition of a 
genus one model. Notice that we modify the definition in the
case $n=2$ to accommodate fields of characteristic $2$.

In \S\ref{invars} we study the ring of invariants. 
We show that it is generated by invariants $c_4$ and $c_6$ 
of the expected degrees and that these invariants 
give a formula for the Jacobian. We also show (in all characteristics) 
that a genus one model defines a smooth curve of genus one if
and only if its discriminant $\Delta = (c_4^3-c_6^2)/1728$
is non-zero. The proofs rely on geometric results proved
in \S\ref{sec:geom} and formulae recorded in \S\ref{sec:wmodels}.

In \S\ref{sec:formulae} we recall some classical methods for
computing the invariants in the cases $n=2,3,4$.
These formulae have already been surveyed in \cite{Mc+},
but are included here to demonstrate our preferred choice of
scalings. In the case $n=5$ we have found an algorithm for
evaluating the invariants. % This is presented in \S\ref{sec:evalalg}.
Our algorithm, presented in \S\ref{sec:evalalg}, 
is inspired by the methods of nineteenth century 
invariant theory, in that we approach the invariants through 
the construction of certain covariants. The key step relies on a 
geometric ``accident'' satisfied by the secant variety 
of a genus one normal quintic. 
% Combining with the
% algorithm in \cite{g1pf} we obtain an algorithm for 
% computing the Jacobian of a genus one normal quintic. 
In \S\ref{compgeom} we compare our invariant-theoretic approach
with some other methods for computing the Jacobian of a genus one curve.

Finally in \S\ref{sec:char23} we include a brief note on 
the invariants in characteristics $2$ and $3$. We find
in these cases that the invariants are insufficient to compute
the Jacobian. Instead it should be possible to find a formula
for the Jacobian that works in all characteristics by
modifying the formulae in characteristic $0$.
This has been carried out by Artin, Rodriguez-Villegas and Tate 
\cite{ARVT} in the case $n=3$.

The formulae and algorithms presented 
in \S\S\ref{sec:formulae},\ref{sec:evalalg} have been contributed 
to {\sf MAGMA} \cite[Version 2.13]{magma} by the author. 
% and will be included in the next release (July 2006).

% \[ \begin{array}{cccccrr}
% n & G_n & \dim X_n & \deg c_4 & \deg c_6 & |c_4| & |c_6| \\ \hline
% 2 & \SL_2 & 5 & 2 & 3 & 3 & 5 \\
% 3 & \SL_3 & 10 & 4 & 6 & 25 & 103 \\
% 4 & \SL_2 \times \SL_4 & 20 & 8 & 12 & 4003 & 81547 \\
% 5 & \SL_5 \times \SL_5 & 50 & 20 & 30 & ? & ? 
% \end{array} \]

\section{Geometric Invariants}
\label{geominvar}

Let $C$ be a smooth curve of genus one defined over $K$,
and let $\omega$ be a non-zero regular 1-form on $C$,
also defined over $K$. 
We say that $\omega$ is an invariant differential.
Over $\Kbar$, the pair $(C,\omega)$ may be put in the form
\begin{equation}
\label{weqn}
y^2 + a_1 xy + a_3 y = x^3 + a_2 x^2 + a_4 x + a_6 
\end{equation}
with $\omega = dx / (2 y + a_1 x + a_3)$.

\begin{Definition}
\label{defginv}
The geometric invariants of the pair $(C,\omega)$ are
\begin{equation*}
% \label{cinv1}
\begin{array}{rcl}
c_4 & = & b_2^2-24 b_4  \\
c_6 & = & -b_2^3 + 36 b_2 b_4 - 216 b_6
\end{array} 
\end{equation*}
where $b_2 = a_1^2+4 a_2$, $b_4 = 2 a_4 + a_1 a_3$ and 
$b_6 = a_3^2 + 4 a_6$.
\end{Definition}

It is clear from the formulae in~\cite[Chapter III]{Si1}
that $c_4$ and $c_6$ depend only on the pair $(C,\omega)$ 
and not on the choice of Weierstrass equation~(\ref{weqn}).
We deduce by Galois theory that $c_4, c_6 \in K$. 

\begin{Lemma} 
\label{geomwts}
If $(C,\omega)$ has geometric invariants $c_4$ and
$c_6$, and $\lambda \in K^*$, then $(C,\lambda^{-1} \omega)$
has geometric invariants $\lambda^4 c_4$ and $\lambda^6 c_6$. 
\end{Lemma}
\begin{Proof} This is again clear from \cite[Chapter III]{Si1}.
\end{Proof}

We show in \S\ref{sec:invdiff} that if a genus one model of degree $n$ 
defines a smooth curve of genus one, then it also defines
an invariant differential on the curve. 
This enables us to construct the invariants $c_4,c_6 \in K[X_n]^{G_n}$ 
as the geometric invariants of the generic genus one model of degree $n$.
In particular we treat the cases $n=2,3,4,5$ in a uniform manner, 
and avoid the problem of combinatorial explosion in the case $n=5$.

The geometric invariants give a formula for the Jacobian.

\begin{Proposition} 
\label{ginvjac}
Assume $\charic(K) \not=2,3$. 
If $(C,\omega)$ has geometric invariants 
$c_4$ and $c_6$ then $C$ has Jacobian
\[ y^2 = x^3 - 27 c_4 x - 54 c_6. \]
\end{Proposition}

The proof relies on two easy lemmas.
 
\begin{Lemma}
\label{om1}
Assume $\charic(K) \not= 2,3$. 
Let $E$ be an elliptic curve defined over $K$ 
with invariant differential $\omega$. 
Let $\alpha$ be an automorphism of $E$. 
Then $\alpha$ is a translation map 
if and only if $\alpha^* \omega = \omega$. 
\end{Lemma}
\begin{Proof}
We write $\tau_P : E \to E$ for translation by $P \in E$. The
map $P \mapsto \tau_P^*(\omega)/\omega$ is a morphism $E \to \Gm$.
It must therefore be constant. Specialising to $P=0$ we deduce
that $\tau_P^*(\omega) = \omega$ for all $P \in E$. (An alternative
proof is given by writing $\tau_P$ as the 
commutator of $\tau_Q$ and $[-1]$ where $2Q=P$.)

Conversely if $\alpha$ is not a translation then $\alpha-1$ is not
constant, and therefore surjective. So $\alpha$ has
a fixed point. Conjugating by a translation, we may suppose that 
the fixed point is $0 \in E$. Since $\charic(K) \not=2,3$ we can put
$E$ in shorter Weierstrass form $y^2 = x^3 + A x + B$. Then the
only automorphisms of $(E,0)$ are of the form $(x,y) \mapsto
(u^2x,u^3y)$. Since $\omega$ is a multiple of $dx/y$ the result is
now clear.
\end{Proof}

\begin{Lemma}
\label{om2}
Assume $\charic(K) \not=2,3$.
Let $E$ be an elliptic curve and $C$ a smooth curve of genus one,
both defined over $K$. 
Let  $\omega_E$ and $\omega_C$ be 
invariant differentials on $E$ and $C$, also defined over $K$.
If there is an isomorphism $\alpha : C \isom E$ defined over $\Kbar$ 
with $\alpha^* \omega_E = \omega_C$ then $E$ is 
the Jacobian of~$C$. 
\end{Lemma}
\begin{Proof}
Let $\xi_\sigma = \sigma(\alpha) \alpha^{-1}$ for 
$\sigma \in \Gal(\Kbar/K)$. Since $\omega_E$ and $\omega_C$ are both
$K$-rational we deduce that $\xi_\sigma^* \omega_E = \omega_E$.
It follows by Lemma~\ref{om1} that $\xi_\sigma$ is a translation.
So $C$ is the twist of $E$ by the class of $\{ \xi_\sigma \}$
in $H^1(K,E)$. In particular $C$ is a torsor under $E$, 
the action $\mu : E \times C \to C$ being given by
\[ \mu(P,Q)  = \alpha(P+ \alpha^{-1}Q). \] 
It follows that $E$ is the Jacobian of $C$.
% Finally $P \mapsto [\mu(P,Q)-Q]$ is an isomorphism $E \isom \Jac(C).$
\end{Proof}

\begin{ProofOf}{Proposition~\ref{ginvjac}}
We are given $(C,\omega)$ with geometric invariants $c_4$
and $c_6$. Let $E$ be the elliptic curve over $K$ with
Weierstrass equation
\[ y^2 = x^3 - 27 c_4 x - 54 c_6. \]
The pairs $(C,\omega)$ and $(E,3 dx/y)$ have the same geometric
invariants, and are therefore isomorphic over $\Kbar$. 
It follows by Lemma~\ref{om2} that $E$ is the Jacobian of $C$. 
\end{ProofOf}

Before we can use Proposition~\ref{ginvjac} to compute the
Jacobian of a genus one curve, we first need to 
compute an invariant differential on the curve.
It is easy to generalise the construction of \S\ref{sec:invdiff} 
to genus one normal curves of arbitrary degree.
An alternative is the following.

Let $C \subset \PP^{n-1}$ be a genus one normal curve of
degree $n$ with hyperplane section $H$. 
We identify the Riemann-Roch space $\CL(H)$ 
with the space of linear forms on $\PP^{n-1}$. If we fix 
$\omega$ then there is a linear map
\[ \begin{array}{cc}
 \wedge^2 \CL(H) \to \CL(2H) \, ; & f \wedge 
g \mapsto \frac{f dg - g df}{\omega}. \end{array} \] 
By Lemma~\ref{projnorm} with $d=2$ 
the natural map $S^2 \CL(H) \to \CL(2H)$ is surjective.
Thus there is an alternating matrix of quadrics % atic forms 
$\Omega = ( \Omega_{ij} )$ with 
\[ \begin{array}{c} 
\omega = \frac{x_j^2 d (x_i/x_j)}{\Omega_{ij}} 
\end{array} \]
for all $ i \not= j$. This matrix has the property that
\[ \begin{pmatrix} \frac{\partial f}{\partial x_1} & \ldots &
\frac{\partial f}{\partial x_n} \end{pmatrix} \Omega \equiv 0 \pmod{I(C)} \]
for all $f \in I(C)$. Starting from generators for $I(C)$
we can use this property to solve for $\Omega$ by linear algebra.
Then $\Omega$ is the data we use to specify $\omega$.
Notice that the entries of $\Omega$ are determined only up to 
the addition of quadrics in $I(C)$. 

\medskip

In \S\ref{compgeom} we compare our invariant-theoretic approach
with some other methods for computing the geometric invariants.

\section{Genus one models}
\label{sec:g1models}

Let $C$ be a smooth curve of genus one defined over $K$, and let 
$D$ be a $K$-rational divisor on $C$ of degree $n$. 
% We go through the 
In each of the cases $n=1,2,3,4,5$
we find equations for the pair $(C,D)$, and use the form of these
equations to motivate our definition of a genus one model.
%  of degree $n$.

\subsection{Genus one models of degree 1}

If $n=1$ then we pick $x,y \in K(C)$ such that $\CL(2D)$ 
and $\CL(3D)$ have bases $1,x$ and $1,x,y$. 
The 7 elements $1,x,y,x^2,xy,x^3,y^2$
in the 6-dimensional space $\CL(6D)$ satisfy a linear dependence relation.
Moreover the coefficients of $x^3$ and $y^2$ are non-zero. Rescaling 
$x$ and $y$ if necessary we find that $C$ has Weierstrass equation 
\begin{equation}
\label{deg1-weqn}
y^2 + a_1 xy + a_3 y = x^3 + a_2 x^2 + a_4 x + a_6. 
\end{equation}

A genus one model of degree $1$ is a tuple 
$\phi = (a_1,a_2,a_3,a_4,a_6)$.
We write $C_\phi \subset \PP^2$ for the curve with Weierstrass 
equation~(\ref{deg1-weqn}). 
% \[ y^2 z + a_1 xyz + a_3 y z^2 = x^3 + a_2 x^2 z + a_4 x z^2 + a_6 z^3. \]
Genus one models $\phi$ and $\phi'$ of degree $1$ are equivalent 
if they are related by substitutions
\[ \begin{array}{rcl}
% z & = & z' \\
x & = & u^2 x' + r \\
y & = & u^3 y' + u^2 s x' + t
\end{array} \]
and $\phi' = u^{-6} \phi$, with $u \not= 0$. 
We write $\G_1$ for the group of all such transformations $[u;r,s,t]$.

\subsection{Genus one models of degree 2}
If $n =2$ then we pick $x,y \in K(C)$ such that $\CL(D)$ 
and $\CL(2D)$ have bases $1,x$ and $1,x,y,x^2$. 
The 9 elements $1,x,x^2,y,x^3,xy,x^4,
x^2 y, y^2$ in the 8-dimensional space $\CL(4D)$ satisfy a 
linear dependence relation. Moreover the coefficient of $y^2$ is non-zero. 
We find that $C$ has equation
\[ y^2 + (\alpha_0 x^2 + \alpha_1 x + \alpha_2 )y 
= a x^4 + b x^3 + c x^2 + d x + e. \]

A genus one model of degree $2$ is a pair of homogeneous 
polynomials $\phi=(p(x,z),q(x,z))$ of degrees 2 and 4.
We write $C_\phi \subset \PP(1,1,2)$ for the curve defined by 
\[ y^2 + p(x,z) y = q(x,z). \]
Here the ambient space $\PP(1,1,2)$ is a weighted project plane,
with degrees $1,1,2$ assigned to the co-ordinates $x,z,y$.
Genus one models $\phi$ and $\phi'$ of degree $2$ are equivalent 
if they are related by substitutions
\[ \begin{array}{rcl}
x & = & B_{11} x' + B_{21} z' \\
z & = & B_{12} x' + B_{22} z' \\
y & = & \mu^{-1} y' + r_0 x'^2 + r_1 x'z' + r_2 z'^2 
\end{array} \]
and $\phi' = \mu^2 \phi$, with $\mu \det B \not= 0$.
We write $\G_2$ for the group of all such transformations 
$[\mu,r,B]$. 

If $\charic(K) \not= 2$ then by completing the square it suffices
to consider models of the form $(0,q(x,z))$. 
These are the binary quartics of Definition~\ref{def:1.1}.

\medskip

If $n \ge 3$ then the complete linear system $|D|$ determines
an embedding $C \to \PP^{n-1}$. We identify $C$ with its image,
which is called a genus one normal curve of degree $n$.
Some basic facts about these curves are recalled in \S\ref{sec:g1nc}. 

\subsection{Genus one models of degree 3}
If $n=3$ then $C \subset \PP^2$ is a plane cubic.
A genus one model of degree 3 is a single homogeneous polynomial
$\phi = (f(x_1,x_2,x_3))$ % with $f$ is a homogeneous polynomial
of degree 3.
% to be a single homogeneous 
% polynomial $\phi = (f(x_1,x_2,x_3))$ with $\deg f = 3$, and
We write $C_\phi \subset \PP^2$ for the variety defined by $f=0$.
Genus one models $\phi$ and $\phi'$ of degree~$3$ are equivalent if they are
related by substitutions $\phi'  = \mu \phi$ and 
$x_j = \sum_{i=1}^3 B_{ij} x_i'$ with $\mu \det B \not= 0$.
We write $\G_3 = \Gm \times \GL_3$ for the group of all such 
transformations. 

\subsection{Genus one models of degree 4}
If $n=4$ then $C \subset \PP^3$ is the complete intersection 
of two quadrics. A genus one model of degree 4 is 
a pair of homogeneous polynomials
\[ \phi = \begin{pmatrix} q_1(x_1,x_2,x_3,x_4) \\ 
 q_2(x_1,x_2,x_3,x_4) \end{pmatrix} \]
of degree $2$. We write $C_\phi \subset \PP^3$ for
the variety defined by $q_1=q_2=0$. 
Genus one models $\phi$ and $\phi'$ of degree~$4$ are equivalent 
if they are related by substitutions $\phi' = A \phi$ and 
$x_j = \sum_{i=1}^4 B_{ij} x_i'$ with $\det A \det B \not= 0$. 
We write $\G_4= \GL_2 \times \GL_4$ for the group of all such 
transformations. 

\subsection{Genus one models of degree 5}
If $n=5$ then $C \subset \PP^4$ is defined by the 
$4 \times 4$ Pfaffians of a $5 \times 5$ alternating matrix
of linear forms.
(See for example \cite{g1pf} and the references cited there.)
A genus one model of degree 5 is a $5 \times 5$ alternating
matrix of linear forms in $5$ variables. We write $C_\phi \subset \PP^4$
for the variety defined by its $4 \times 4$ Pfaffians.
Genus one models $\phi$ and $\phi'$ of degree~$5$ are equivalent if they are
related by substitutions $\phi' = A \phi A^T$ 
and $x_j = \sum_{i=1}^4 B_{ij} x_i'$ with $\det A \det B \not= 0$. 
We write $\G_5= \GL_5 \times \GL_5$ for the group of all such 
transformations. 

\section{The ring of invariants}
\label{invars}

Let $X_n$ be the space of genus one models of degree $n$.
For $n=1,2,3,4,5$ this is an affine space of
dimension $N= 5,8,10,20,50$. The co-ordinate ring $K[X_n]$ 
is a polynomial ring in $N$ variables. 
For $n=3,4,5$ we give this ring its usual grading by degree.
In the cases $n=1,2$ the rings are
\[ \begin{array}{rcl}
K[X_1] & = & K[a_1, a_2, a_3, a_4, a_6] \\
K[X_2] & = & K[\alpha_0, \alpha_1, \alpha_2, a,b,c,d,e]. \\
\end{array} \]
We assign degrees $\deg(a_i)=i$, $\deg(\alpha_i)=1$, 
$\deg(a)= \ldots = \deg(e)=2$.
In \S\ref{sec:g1models} we defined a linear algebraic group 
$\G_n$ acting on $X_n$. We now write $G_n$ for the commutator subgroup,
{\em i.e.}
\[ \begin{array}{rcl}
G_1 & = & \{ [1;r,s,t] \in \G_1 \} \\ % : r,s,t \in \Ga \} \\
G_2 & = & \{ [1,r,B] \in \G_2 :  B \in \SL_2 \} \\ % r \in \Ga^3,
G_3 & = & \SL_3 \\
G_4 & = & \SL_2 \times \SL_4 \\
G_5 & = & \SL_5 \times \SL_5.
\end{array} \]
 
\begin{Definition}
The ring of invariants is 
\[ K[X_n]^{G_n} = \{ F \in K[X_n] : F \circ g = F 
\text{ for all } g \in G_n(\Kbar) \}. \]
The definition is extended to an integral domain $R$ by putting
% If $R$ is an integral domain with field of fractions $K$ then 
\[R[X_n]^{G_n} = R[X_n] \cap K[X_n]^{G_n}.\]
where $K$ is the field of fractions of $R$.
\end{Definition}

We define a rational character $\det : \G_n \to \Gm$
\[ \begin{array}{l@{\qquad}rcl}
n = 1 & [u;r,s,t] & \mapsto & u^{-1} \\
n = 2 & [\mu,r,B] & \mapsto & \mu \det B \\ 
n = 3 & [\mu,B] & \mapsto & \mu \det B \\ 
n = 4 & [A,B] & \mapsto & \det A \det B \\
n = 5 & [A,B] & \mapsto & (\det A)^2 \det B. 
\end{array} \]
\begin{Definition}
The vector space of invariants of weight $k$ is
\[ K[X_n]_k^{G_n} = \{ F \in K[X_n] : F \circ g = (\det g)^k F 
\text{ for all } g \in \G_n(\Kbar) \}. \]
\end{Definition}

\begin{Lemma} 
\label{degwt}
Every % If $F \in K[X_n]^{G_n}$ is a 
homogeneous invariant of degree~$d$ is an invariant of weight $k$ where
\[ d=  \left\{ \begin{array}{ll} k & \text{ if } n =1,2,3 \\
2k & \text{ if } n=4 \\ 5k & \text{ if } n=5. 
\end{array} \right.  \]
In particular the ring of invariants is graded by weight, {\em i.e.}
\[K[X_n]^{G_n} = \oplus_{k \ge 0} K[X_n]_k^{G_n}.\]
\end{Lemma}
\begin{Proof} 
We treat the cases $n=4,5$.
% The inclusion ``$\supset$'' in (\ref{grading}) is clear.
% We prove ``$\subset$'' in the cases $n=4,5$.
% Let $F \in K[X_n]^{G_n}$ be a homogeneous invariant of degree $d$.
 Since the only rational characters of $\G_n$
are of the form $[A,B] \mapsto (\det A)^p (\det B)^q$ we have
 \[  F \circ [A,B] =  (\det A)^p (\det B)^q  F  \]
for some integers $p,q$.
Considering $[A,B]$ in the centre of $\G_n$ we deduce 
\[ \begin{array}{lrcl} \smallskip
n=4 & \left\{ \begin{array}{rcl}
d & = & 2p  \\
2 d  & = & 4q  
\end{array} \right.  \\
n=5 & \left\{ \begin{array}{rcl}
2 d & = & 5p \\
 d & = & 5q 
\end{array} \right.  
\end{array} \]
We are done by the definition of $\det : \G_n \to \Gm$. 
The cases $n=1,2,3$ are similar.
\end{Proof}

We are ready to state our main theorem.
\begin{Theorem}
\label{mainthm}
There are invariants $c_4, c_6, \Delta \in K[X_n]^{G_n}$ 
of weights $4$, $6$ and $12$, related by $c_4^3 - c_6^2 = 1728 \Delta$, 
such that \\
(i) If $\charic(K) \not= 2,3$ then $K[X_n]^{G_n} = K[c_4,c_6]$. \\
(ii) The variety $C_\phi$ defined by $\phi \in X_n$ is a smooth curve of
genus one if and only if $\Delta(\phi)\not=0$. \\
(iii) If $\charic(K) \not= 2,3$ and $\phi \in X_n$ with $\Delta(\phi) \not=0$ 
then $C_\phi$ has Jacobian
\[ y^2 = x^3 - 27 c_4(\phi) x - 54 c_6(\phi). \]
\end{Theorem}

The proof depends on the following geometric statements.

\begin{Proposition}
\label{PropA}
Assume $K=\Kbar$. Let $X_n^{\sing}$ be the set of all 
models $\phi \in X_n$ which do not define a smooth curve of genus one. 
% Let $X_n^{\sing}$ be the set of all $\phi \in X_n$ for
% which $C_\phi$ is not a smooth curve of genus one. 
Then $X_n^{\sing}$ is an irreducible Zariski closed subset of $X_n$. 
In particular the generic genus one model of degree $n$ 
defines a smooth curve of genus one.
\end{Proposition}
\begin{Proof}
The cases $n=1,2,3$ are well known.
A proof for  $n=3,4,5$ is given in \S\ref{sec:gen}.
\end{Proof}

Let $\PP(X_n)$ be the projective space determined by $X_n$. (This is a
weighted projective space in the cases $n=1,2$.) We recall that 
elements of $X_n$ are equivalent if they lie in the same $\G_n$-orbit. 

\begin{Proposition}
\label{PropB}
Assume $K=\Kbar$. Let $\phi,\phi' \in X_n$ with $C_\phi$ and $C_{\phi'}$
either smooth curves of genus one or rational curves with a single node.
Then $C_\phi$ and $C_{\phi'}$ are isomorphic as curves if and only if 
$\phi$ and $\phi'$ are equivalent. Moreover the stabiliser of $\phi$
for the action of $G_n$ on $\PP(X_n)$ % $\{g \in G_n :g \phi = \phi \}$ 
is finite. 
\end{Proposition}
\begin{Proof}
The cases $n=1,2$ are straightforward.
A proof for  $n=3,4,5$ is given in \S\ref{minfreeres}.
\end{Proof}

We identify $K[X_n]^{G_n}$ as a subring of $K[X_1]^{G_1}$. To do this
we start with an elliptic curve $E$ in Weierstrass form
\[ y^2 + a_1 x y + a_3 y = x^3 + a_2 x^2 + a_4 x + a_6. \]
The complete linear system $|n.0|$ determines a morphism
$E \to \PP^{n-1}$. 
% \[ \begin{array}{llcl}
% n=2 \qquad  & (x_0:x_2) & = & (1:x) \\
% n=3 & (x_0:x_2:x_3) & = & (1:x:y) \\
% n=4 & (x_0:x_2:x_3:x_4) & = & (1:x:y:x^2) \\
% n=5 & (x_0:x_2:x_3:x_4:x_5) & = & (1:x:y:x^2:xy).
% \end{array} \]
The image is described by a genus one model of degree $n$. 
In \S\ref{sec:wmodels} we specify such a model and hence define
a morphism $\pi_n: X_1 \to X_n$. The models $\pi_n(\phi)$ 
for $\phi \in X_1$ are called Weierstrass models.
Collectively they form the Weierstrass family.

\begin{Proposition}
\label{PropC}
There are morphisms $\pi_n: X_1 \to X_n$ and $\gamma_n : \G_1 \to \G_n$
with the following properties. \\
(i) If $\phi' = \pi_n(\phi)$ then $C_\phi$ and $C_{\phi'}$ are isomorphic
as curves. \\
(ii) $\gamma_n$ is a group homomorphism. \\
(iii) $(\gamma_n g)(\pi_n \phi) = \pi_n(g \phi)$ for all $g \in \G_1$ and
$\phi \in X_1$. \\
(iv) $\det(\gamma_n g)= \det(g)$ for all $g \in \G_1$. 
\end{Proposition}
\begin{Proof}
The proposition is checked by direct computation using 
the formulae in \S\ref{sec:wmodels}. 
\end{Proof}

The map $\pi_n : X_1 \to X_n$ induces a homomorphism
of polynomial rings $\pi_n^*: K[X_n] \to K[X_1] \, ; \,
 F \mapsto F \circ \pi_n$. By Proposition~\ref{PropC} 
it restricts to a homomorphism of graded rings
\[ \pi_n^* : K[X_n]^{G_n} \to K[X_1]^{G_1} \]
where the grading is by weight.

\begin{Lemma} 
\label{INJ}
The map $\pi_n^* : K[X_n]^{G_n} \to K[X_1]^{G_1}$
is an injective homomorphism of graded rings.
\end{Lemma}
\begin{Proof} 
We must show that $\pi_n^*$ is injective.
For this we are free to assume that $K$ is algebraically closed.
If $F \in K[X_n]^{G_n}$ is a homogeneous invariant vanishing on the 
Weierstrass family then by Propositions~\ref{PropB} and~\ref{PropC} 
it also vanishes at every $\phi \in X_n$ for which $C_\phi$ 
is a smooth curve of genus one.
% \[ \{ \phi \in X_n : C_\phi \text{ is a smooth curve of genus one }\}. \]
Proposition~\ref{PropA} tells us that the latter are Zariski dense
in $X_n$. It follows that $F$ is identically zero and hence 
$\pi_n^*$ is injective.
\end{Proof}

% We show that $\pi_n^*$ is an isomorphism, although the proof
% in the case $\charic(K) =2,3$ is postponed to \S\ref{charK23}.
% First we find generators for 
Computing the ring $K[X_1]^{G_1}$ is entirely routine.
We recall that \[K[X_1]= K[a_1,a_2,a_3,a_4,a_6].\]
Following Tate's formulaire \cite[Chapter III]{Si1} we put
\begin{equation}
\label{binv}
\begin{array}{rcl}
b_2 & = & a_1^2+4 a_2 \\
b_4 & = & 2 a_4 + a_1 a_3 \\
b_6 & = & a_3^2 + 4 a_6 \\
b_8 & = & a_1^2 a_6 + 4 a_2 a_6 - a_1 a_3 a_4 + a_2 a_3^2 - a_4^2
\end{array} 
\end{equation}
and 
\begin{equation}
\label{cinv}
\begin{array}{rcl}
c_4 & = & b_2^2-24 b_4  \\
c_6 & = & -b_2^3 + 36 b_2 b_4 - 216 b_6 \\
\Delta & = & -b_2^2 b_8 - 8 b_4^3 - 27 b_6^2 + 9 b_2 b_4 b_6. 
\end{array} 
\end{equation}
It is well known that $c_4,c_6,\Delta \in \Z[X_1]^{G_1}$
and $c_4^3-c_6^2 = 1728 \Delta$.
\begin{Lemma} 
\label{inv1}
If $\charic(K) \not= 2,3$ then $K[X_1]^{G_1}
=K[c_4,c_6]$.
\end{Lemma}
\begin{Proof}
This is Theorem~\ref{mainthm}(i) in the case $n=1$.
It is an immediate consequence of the standard procedure
for putting a Weierstrass equation in the shorter form 
$y^2= x^3 + Ax + B$. The required isomorphism is $\iota^*$ where 
\[ \begin{array}{cc} \iota : \Aff^2 \to X_1 \, ; & (c_4,c_6) 
 \mapsto  (0,0,0,-c_4/ 48 ,-c_6/864). \end{array} \]
\end{Proof}

We have reduced the proof of Theorem~\ref{mainthm}(i) to showing 
that $\pi_n^*$ is surjective. Equivalently, we must show that 
$K[X_n]^{G_n}$ contains invariants of weights 4 and 6.
% now reduces to showing that
% $\pi_n^*$ is surjective. In other words we must construct some invariants.
One method would be to split into the cases $n=2,3,4,5$ and use
the explicit constructions 
presented in \S\S\ref{sec:formulae},\ref{sec:evalalg}. 
This makes the theorem appear an accident, especially in the case 
$n=5$. Instead we give a construction based on Proposition~\ref{PropB}.

\begin{Lemma} 
\label{ORB}
Assume $K=\Kbar$. Let $\phi,\phi' \in X_n$ with $C_{\phi}$ and $C_{\phi'}$
either smooth curves of genus one or rational curves with a single node.
Then the Zariski closure of the $\G_n$-orbit of $\phi$ is the zero locus
of an irreducible homogeneous invariant $F \in K[X_n]^{G_n}$.
Moreover $F(\phi')=0$ if and only if $\phi$ and $\phi'$ are equivalent.
\end{Lemma}
\begin{Proof}
By Proposition~\ref{PropB} the morphism
$G_n \to \PP(X_n) ; \, g \mapsto g (\phi)$
has zero-dimensional fibres. 
% \[ \begin{array}{ccccc}
% n & & G_n &  \dim G_n & \dim X_n \\ \hline
% 2 & & \SL_2 \ltimes \Ga^3 &  6 & 8   \\
% 3 & & \SL_3               &  8 & 10  \\
% 4 & & \SL_2 \times \SL_4  & 18 & 20  \\
% 5 & & \SL_5 \times \SL_5  & 48 & 50 
% \end{array} \]
% \end{Proof}
But for each $n$ we find \[\dim(G_n) = \dim(X_n)-2.\]
So the $G_n$-orbit of $\phi$ in $\PP(X_n)$ has codimension 1. Moreover
since $G_n$ is irreducible, every $G_n$-orbit is irreducible. Therefore 
the Zariski closure of the orbit of $\phi$ is the zero locus of an
irreducible homogeneous polynomial $F \in K[X_n]$. Since the equivalence
class of $\phi$ determines $F$ uniquely up to scalars, 
and $G_n$ is the commutator subgroup of $\G_n$, 
it follows that $F$ is an invariant.

If $\phi$ and $\phi'$ are equivalent then clearly $F(\phi')=0$.
For the converse we suppose $F(\phi')=0$. Then the $G_n$-orbits of 
$\phi$ and $\phi'$ in $\PP(X_n)$ have the same Zariski closure, $Z$ say.
A standard argument (see e.g. \cite[Chapter I, \S5.3, Theorem 6]{Sh}) shows
that each of these orbits contains a non-empty open subset of $Z$. 
Since $Z$ is irreducible these open sets must intersect. It follows that
$\phi$ and $\phi'$ are equivalent.
\end{Proof}

We restrict these invariants to the Weierstrass family.

\begin{Lemma}
\label{POW}
Assume $K= \Kbar$ and $\charic(K) \not= 2,3$.
Then there are irreducible invariants $F_4,F_6 \in K[X_n]^{G_n}$
and integers $p,q \ge 1$ such that $F_4 \circ \pi_n = c_4^p$ 
and $F_6 \circ \pi_n = c_6^q$.
\end{Lemma}
\begin{Proof}
By Proposition~\ref{PropC} we can pick $\phi \in X_n$ with $C_\phi$ a
smooth curve of genus one with $j$-invariant $0$.
Let $F_4 \in K[X_n]^{G_n}$ be the invariant constructed from $\phi$
in Lemma~\ref{ORB}. Then $F_4 \circ \pi_n$ is a homogeneous element
of $K[X_1]^{G_1} = K[c_4,c_6]$. Rescaling $F_4$ we can write
\[ F_4 \circ \pi_n = c_4^p \,\, c_6^q \,\, \Delta^r \, 
\prod_{\nu=1}^s (c_4^3-j_\nu \Delta) \]
for some integers $p,q,r,s \ge 0$ and constants $j_1, \ldots j_s \not=0,1728$.

Now let $\phi' = \pi_n(\phi_1)$ be a Weierstrass model 
with $C_{\phi'}$ a smooth curve of genus one.
If this curve has $j$-invariant not equal to $0$ then by
Lemma~\ref{ORB} we have $F_4(\phi') \not=0$.
By varying the choice of $\phi_1$ we deduce that $q=s=0$. 
We then repeat the argument for $C_{\phi'}$ a Weierstrass model 
with a node. This shows that $r=0$.
% or a rational curve with a node. Lemma~\ref{ORB} tells us that 
% $F_4(\phi') \not=0$ and from this we deduce $q=r=s=0$. 
% Hence $F_4 \circ \pi_n = c_4^p$. 
The statement for $c_6$ is proved similarly, 
starting with $j$-invariant 1728.
\end{Proof}

The proof of Theorem~\ref{mainthm}(i) now reduces to showing
that $p=q=1$ in Lemma~\ref{POW}. For this we quote 
a geometric result whose proof uses properties of the
invariant differential. 

\begin{Definition}
\label{def:propeq}
Genus one models $\phi, \phi' \in X_n$ are properly equivalent 
if there exists $g \in \G_n$ with $g \phi = \phi'$ and $\det(g)=1$.
\end{Definition}

\begin{Proposition}
\label{PropD}
Assume $K=\Kbar$ and $\charic(K) \not= 2,3$. 
Let $\phi \in X_n$ with $C_\phi$ a smooth curve of genus one. 
Then $\phi$ is properly equivalent to $\pi_n(0,0,0,A,B)$ for 
some unique $A,B \in K$. 
\end{Proposition}
\begin{Proof} The existence is already clear from Propositions~\ref{PropB} 
and~\ref{PropC}. We prove uniqueness in \S\ref{sec:invdiff}.
\end{Proof}

\begin{Lemma}
\label{KK}
Assume $K= \Kbar$ and $\charic(K) =0$. Then the map $\pi_n^* : K[X_n]^{G_n}
\to K[X_1]^{G_1}$ is surjective.
\end{Lemma}
\begin{Proof}
%  By Lemma~\ref{POW} there are irreducible invariants 
%  $F_4,F_6 \in K[X_n]^{G_n}$ and integers $p,q \ge 1$ with 
%  $F_4 \circ \pi_n = c_4^p$ and $F_6 \circ \pi_n = c_6^q$.
Let $\phi \in X_n(\K)$ be the generic model defined over the function 
field $\K = K(X_n)$. We have assumed $\charic(K)=0$ so that $\K$ is perfect. 
Proposition~\ref{PropA} tells us that $C_\phi$ is a smooth curve of
genus one. So by Proposition~\ref{PropD}, 
$\phi$ is properly equivalent to $\pi_n(0,0,0,A,B)$ for some unique 
$A,B \in \overline{\K}$. The uniquess statement shows that $A$ and $B$
are fixed by $\Gal(\overline{\K}/\K)$ and hence $A,B \in \K$. 

Let $F_4,F_6 \in K[X_n]^{G_n}$ be the  irreducible invariants 
constructed in Lemma~\ref{POW}.
Then $F_4 = f_4^p$ and $F_6= f_6^q$ where
\[ \begin{array}{rcccl}
f_4 & = & c_4(0,0,0,A,B) & = & -48 A \\
f_6 & = & c_6(0,0,0,A,B) & = & -864 B.
\end{array} \]
Since $K[X_n]$ is integrally closed (in its field of fractions $\K$) 
and $F_4$, $F_6 \in K[X_n]$ are irreducible it follows that $p=q=1$.
\end{Proof}

Applying Lemma~\ref{KK} with $K= \Qbar$ we learn that the invariants 
$c_4,c_6,\Delta \in \Z[X_1]^{G_1}$ extend to invariants in $\Qbar[X_n]^{G_n}$.
These invariants are again denoted $c_4,c_6,\Delta$.
Since $\pi_n^*$ is injective it follows by Galois theory that 
$c_4,c_6,\Delta  \in \Q[X_n]$. In fact the coefficients are integers.
\begin{Lemma} 
\label{intcoeff}
$c_4,c_6, \Delta \in \Z[X_n]$.
\end{Lemma}
\begin{Proof}
Let $F=c_4, c_6$ or $\Delta$. Let $p$ be a prime and $r \ge 0$ 
an integer. We suppose for a contradiction that $p^{r+1}F 
\in \Z_p[X_n]$ yet $p^r F \not\in \Z_p[X_n]$. Then each
coefficient of $\pi_n^*(p^{r+1}F) \in \Z[X_1]$ is divisible by $p$.
So if $G \in \F_p[X_n]$ is the reduction of 
$p^{r+1}F$ mod $p$ then $\pi_n^* G = 0$. The injectivity 
established in Lemma~\ref{INJ} shows that $G=0$. 
Therefore $p^r F \in \Z_p[X_n]$. This is the required
contradiction.
\end{Proof}

\begin{Remark} 
\label{primrem}
Since the original $c_4, c_6, \Delta \in \Z[X_1]$ are primitive
it is clear that the new $c_4, c_6, \Delta \in \Z[X_n]$ are
also primitive. This means it is possible to specify our scalings
of $c_4, c_6, \Delta$, at least up to sign, without the need to
compute their restrictions to the Weierstrass family.
\end{Remark}

We revert to working over an arbitrary perfect field $K$.

\medskip

\begin{ProofOf}{Theorem~\ref{mainthm}}
Let $c_4,c_6, \Delta \in K[X_n]$ be the images of 
$c_4,c_6,\Delta \in \Z[X_n]$. 
These polynomials are invariants of weights $4$, $6$ and $12$,
satisfying $c_4^3- c_6^2 = 1728 \Delta$.
They are non-zero by Remark~\ref{primrem}. \\
(i) If $\charic(K) \not= 2,3$ then by Lemmas~\ref{INJ} and~\ref{inv1} 
the map \[\pi_n^* : K[X_n]^{G_n} \to K[X_1]^{G_1}  = K[c_4,c_6]\]
is an isomorphism of graded rings. \\ 
% In particular  $K[X_n]^{G_n} = K[c_4,c_6]$. \\
(ii) We may assume that $K$ is algebraically closed.
If $\phi \in X_n$ with $C_\phi$ a smooth curve of genus one 
then by Propositions~\ref{PropB} and~\ref{PropC} 
it is equivalent to a Weierstrass model.
% to find $\phi_1 \in X_1$ with $\pi_n(\phi_1)$ equivalent to $\phi$.
We deduce $\Delta(\phi) \not= 0$.
So there is an inclusion 
\[ \{ \Delta=0 \} \subset  X_n^{\sing}. \]
But Proposition~\ref{PropA} asserts that $X_n^{\sing}$ 
is closed and irreducible. So the inclusion is in fact an equality. \\
% We complete the proof of Theorem~\ref{mainthm}.
(iii) Let $\phi \in X_n$ with $C_\phi$ a smooth curve of genus one.
In \S\ref{sec:invdiff} we use $\phi \in X_n$ to define
an invariant differential $\omega_\phi$ on $C_\phi$. 
In the case $n=5$ we must assume $\charic(K) \not=2$. We then show
that $c_4(\phi)$ and $c_6(\phi)$ are the geometric invariants
of the pair $(C_\phi,\omega_\phi)$. The formula for the
Jacobian follows by Proposition~\ref{ginvjac}.
\end{ProofOf}

Theorem~\ref{mainthm}(iii) is proved in \cite{Mc+} for $n=2,3,4$ 
by giving formulae for the covering map (of degree $n^2$) from 
a genus one curve to its Jacobian. We have extended 
to the case $n=5$ by taking a different approach, based on 
properties of the invariant differential. 

It turns out that the map $\pi_n^* : K[X_n]^{G_n} \to K[X_1]^{G_1}$ 
is an isomorphism in all characteristics. The proof in 
characteristics 2 and 3 is given in \S\ref{sec:char23}.

\medskip

The remaining sections of the paper may be read in any order.

\section{Geometry}
\label{sec:geom}

The aim of this section is to prove the geometric results 
cited in \S\ref{invars}. We work over an algebraically closed field $K$.
The homogeneous ideal of a projective variety $X$ is denoted $I(X)$.

\subsection{Genus one normal curves}
\label{sec:g1nc}
We recall some basic facts about genus one normal curves
and rational nodal curves.

% We are interested in projective curves of the following two kinds.
% the following projective curves.

\begin{Definition}
\label{def:g1ncrnc}
Let $n \ge 3$ be an integer. \\
(i) A genus one normal curve $C \subset \PP^{n-1}$ is a smooth
curve of genus one and degree $n$ that spans $\PP^{n-1}$. \\
(ii) A rational nodal curve $C \subset \PP^{n-1}$ is a rational
curve of degree $n$ that spans $\PP^{n-1}$ and has a single node. 
\end{Definition}

\begin{Remark}
\label{rem:linsys}
Equivalently, a genus one normal curve is a smooth 
curve of genus one embedded by a complete linear
system of degree~$n$. A rational nodal curve is the 
image of a morphism $\PP^1 \to \PP^{n-1}$ determined by a linear 
system of the form
\[ \{ f \in \CL(D) | f(P_1)=f(P_2) \} \]
for $D$ a divisor on $\PP^1$ of degree $n$, and $P_1$, $P_2 \in \PP^1$ 
distinct. 
\end{Remark}

\begin{Proposition}
\label{nquads}
Let $C \subset \PP^{n-1}$ be either a genus one normal curve
or a rational nodal curve. 
If $n \ge 4$ then the ideal $I(C)$ is generated 
by a vector space of quadrics of dimension $n(n-3)/2$. 
\end{Proposition}

This proposition is well known, 
at least for genus one normal curves. 
Our proof, based on an argument in \cite{Hu}, has the 
advantage of working for rational nodal curves at the same time.
We write $R = K[x_1, \ldots ,x_{n}]$ and $R' = K[x_1, \ldots ,x_{n-1}]$
for the homogeneous co-ordinate rings of $\PP^{n-1}$ and 
$\PP^{n-2}$. We give each ring its usual grading by degree,
say $R = \oplus_{d \ge 0} R_d$ and $R' = \oplus_{d \ge 0} R'_d$.
\begin{Lemma}
\label{lem4}
Let $X \subset \PP^{n-2}$ be a set of $n$ points in general position. \\
(i) The evaluation map $\pi_X: R'_d \to K^n$
is surjective for all $d \ge 2$. \\
(ii) If $n \ge 4$ then the ideal $I(X) \subset R'$ is generated by quadrics. 
\end{Lemma}
\begin{Proof}
We change co-ordinates so that $X$ is the set of points
$(1:0: \ldots :0)$, $(0:1: \ldots :0)$, \ldots, $(0:0: \ldots :1)$
and $(1:1: \ldots :1)$. The proof is now straightforward.
\end{Proof}
We show that the curves defined in Definition~\ref{def:g1ncrnc} 
are projectively normal.
\begin{Lemma}
\label{projnorm}
Let $C \subset \PP^{n-1}$ be either a genus one normal curve
or a rational nodal curve. Let $H$ be the divisor of a hyperplane
section, say cut out by a linear form $h \in R_1$.
%  be a linear form
% cutting out a divisor $H$ on $C$. 
Then the map
\[\pi_{C}: R_d \to \CL(dH) \,; \,\, f \mapsto f/h^d \]  
is surjective for all $d \ge 1$. 
\end{Lemma}
\begin{Proof} The proof is by induction on $d$,
the case $d=1$ being clear from Riemann-Roch.
For the induction step we choose a hyperplane $\{\xi=0\}$ meeting $C$ 
in $n$ distinct points disjoint from $H$.
Again by Riemann-Roch any $n-1$ distinct points
on $C$ span a hyperplane. So $X=C \cap \{\xi=0\}$ satisfies the
hypothesis of Lemma~\ref{lem4}.

Let $d \ge 2$. We are given $f \in \CL(dH)$ 
and wish to show that it belongs to the image of $\pi_{C}$.
By Lemma~\ref{lem4}(i) it
suffices to treat the case where $f$ vanishes on $X$. 
But then $f=(\xi/h)f'$ for some $f' \in \CL((d-1)H)$.
Applying the induction hypothesis to $f'$, we deduce that
$f$ is in the image of $\pi_{C}$ as required.
\end{Proof}

\begin{ProofOf}{Proposition~\ref{nquads}}
% The case $n=3$ being clear, we assume $n \ge 4$.
We continue with the notation of the last proof.
Since $C$ is contained in no hyperplane, the natural map
\begin{equation}
\label{cxisom}
I(C) \cap R_2 \to I(X) \cap R'_2
\end{equation}
is injective. By Lemmas~\ref{lem4} and~\ref{projnorm} these
spaces each have dimension $n(n-3)/2$. So~(\ref{cxisom}) is an
isomorphism.  
Now let $f \in I(C) \cap R_d$. We must show that $f$ is in 
the ideal generated by $I(C) \cap R_2$.
By Lemma~\ref{lem4}(ii) and the surjectivity of (\ref{cxisom}) it suffices 
to treat the case where $f$ vanishes on $X$. But then $f=\xi f'$ for some
$f' \in I(C) \cap R_{d-1}$. The proposition now follows
by induction on $d$.
\end{ProofOf}

We say that curves $C,C' \subset \PP^{n-1}$ are
projectively equivalent if there exists $\alpha \in \PGL_n$
with $\alpha(C)=C'$.

\begin{Lemma}
\label{trans}
(i) Genus one normal curves $C,C' \subset \PP^{n-1}$
are projectively equivalent if and only if they have the
same $j$-invariant. \\
(ii) Any two rational nodal curves $C,C' \subset \PP^{n-1}$ are
projectively equivalent.
\end{Lemma}
\begin{Proof}
(i) Let $C$ and $C'$ have hyperplane sections $H$ and $H'$.
If $C$ and $C'$ are isomorphic as curves then composing with a 
translation map we can find an isomorphism 
$\alpha : C \isom C'$ with $\alpha^* H' \sim H$. \\
(ii) This is clear from Remark~\ref{rem:linsys}.
\end{Proof}

\begin{Lemma}
\label{finite}
Let $C \subset \PP^{n-1}$ be either a genus one normal curve or a rational
nodal curve. Then there are only finitely many $\alpha \in \PGL_n$
with $\alpha(C)=C$. 
\end{Lemma}
\begin{Proof}
We first treat the case $C$ is a genus one normal curve, 
say with hyperplane section $H$. 
We are interested in the automorphisms $\alpha$ of $C$ 
with $\alpha^* H \sim H$. %  where $H$ is the hyperplane section.
% Let $E$ be the Jacobian of $C$.
The automorphism group of $C$ sits in an exact sequence
\[ 0 \to E \to \Aut(C) \to \Aut(E,0) \to 0 \]
where $E$ is the Jacobian of $C$. 
The first map is $P \mapsto \tau_P$ where $\tau_P$ is translation by $P$.
Since $H$ is a divisor of degree $n$ we have $\tau_P^* H \sim H$ 
if and only if $nP =0$. The lemma follows from the fact that
$E[n]$ and $\Aut(E,0)$ are both finite.

If $C$ is a rational nodal curve then without loss of generality
it is the image of
\[ \PP^1 \mapsto \PP^{n-1} \, ; \quad
(s:t) \mapsto (s^n+t^n : s t^{n-1} : \ldots : s^{n-1}t). \]
The group of automorphisms of $\PP^1$ that extend to automorphisms
of $\PP^{n-1}$ form a copy of the dihedral group 
generated by $(s:t) \mapsto (t:s)$ and $(s:t) \mapsto ( \zeta s : t)$ 
for $\zeta$ an $n$th root of unity.
\end{Proof}

\subsection{Minimal free resolutions}
\label{minfreeres}
We recall that a genus one model of degree $n=3,4,5$ is a 
collection of homogenoeus polynomials in $R= K[x_1, \ldots,x_n]$.
Splitting into the cases $n=3,4,5$ we now use
$\phi \in X_n$ to define an ideal $I_\phi \subset R$
and a complex of graded free $R$-modules $\FF_{\bullet}(\phi)$.
We write $R(d)$ for the graded $R$-module with $R(d)_e = R_{d+e}$.

If $n=3$ then $\phi$ consists of a single polynomial $f \in R$.
This polynomial generates an ideal $I_\phi \subset R$ and defines a complex
\[ \FF_{\bullet}(\phi) : \qquad 0 \ra R(-3) \stackrel{f}{\ra} R \ra 0. \]

If $n=4$ then $\phi$ consists of polynomials $q_1,q_2$.
These polynomials generate an ideal $I_\phi \subset R$ and define a complex
\[ \FF_{\bullet}(\phi) : \qquad 
0 \ra R(-4) \stackrel{ \begin{pmatrix} -q_2 \\ q_1
\end{pmatrix}}{\ra} R(-2)^2 \stackrel{ \begin{pmatrix} q_1 & q_2
\end{pmatrix} }{\ra} R \ra 0. \]

If $n=5$ then $\phi$ is a $5 \times 5$ alternating matrix of linear forms.
The Pfaffian of a $4 \times 4$ alternating matrix is
\[ \pf \begin{pmatrix} 0 & a_1 & a_2 & a_3 \\
& 0 & b_3 & b_2 \\  &  & 0 & b_1 \\ & & & 0 \end{pmatrix}
= a_1 b_1 - a_2 b_2 + a_3 b_3. \]
We write $\phi^{\{i\}}$ for the submatrix of $\phi$ obtained 
by deleting the $i$th row and $i$th column.
Then the vector of submaximal Pfaffians of $\phi$ is $P=(p_1, \ldots ,p_5)$
where \[p_i= (-1)^{i+1} \pf(\phi^{\{i\}}).\]
% and $\phi^{\{i\}}$ is the submatrix obtained by deleting the 
% $i$th row and $i$th column of $\phi$.
These polynomials generate an ideal $I_\phi \subset R$ and define a complex
% The quadrics $p_1, \ldots, p_5$ generate 
% an ideal $I_\phi \subset R$. We define a complex
\[ \FF_{\bullet}(\phi) : \quad 
0 \ra R(-5) \stackrel{P^T}{\ra} R(-3)^5 \stackrel{\phi}{\ra} 
R(-2)^5 \stackrel{P}{\ra} R \ra 0. \]

In each case $n=3,4,5$, the variety $C_\phi \subset \PP^{n-1}$ is 
that defined by the ideal $I_\phi \subset R$. We say that
$\FF_{\bullet}(\phi)$ is a minimal free resolution of $R/I_\phi$ 
if it is exact at every term except the final copy of $R$ where 
the homology is $R/I_\phi$.
\begin{Lemma} 
\label{isminres}
Let $n=3,4,5$ and let $\phi \in X_n$.  \\
% Let $\phi \in X_n$ be a genus one model of degree $n=3,4,5$. \\
(i) Every component of $C_\phi$ has dimension at least $1$. \\
(ii) If every component of $C_\phi$ has dimension $1$ then 
$\FF_{\bullet}(\phi)$ is a minimal free resolution of $R/I_\phi$. 
\end{Lemma}
\begin{Proof}
(i) This is clear for $n=3,4$. For $n=5$ we recall that 
the $4 \times 4$ Pfaffians of a generic $5 \times 5$ alternating 
matrix define the image of the Plucker embedding $\Gr(2,5) \to \PP^9$.
Then $C_\phi$ is the intersection of this
Grassmannian with a linear subspace $\PP^4$. Since $\Gr(2,5)$ 
has dimension~6 we are done by \cite[I, Theorem 7.2]{Ha}. \\
(ii) If $n=3,4$ then our claim is that $f$ is non-zero, respectively
that $q_1,q_2$ are coprime. This is clear. The case $n=5$ is
an application of the Buchsbaum-Eisenbud acyclicity criterion,
for which we refer to \cite[Theorem 1.4.13]{CMrings} 
or \cite[Theorem 20.9]{E}.
\end{Proof}

We recall that if $A$ is a finitely generated graded $K$-algebra, 
say $A = \oplus_{d \ge 0} A_d$, then there is a polynomial $h_{A}(t)$,
called the Hilbert polynomial, with the property that $h_{A}(d) = \dim(A_d)$
for all $d \gg 0$.

\begin{Lemma} 
\label{hilbpoly}
(i) Let $n=3,4,5$ and let $\phi \in X_n$. 
If the complex $\FF_{\bullet}(\phi)$ 
is a minimal free resolution of $R/I_\phi$ then \[h_{R/I_\phi}(t) = nt.\] 
(ii) If $C \subset \PP^{n-1}$ is a curve of arithmetic genus $g$ 
and degree $d$ then 
\[ h_{R/I(C)}(t) = dt+(1-g).\]
\end{Lemma}
\begin{Proof}
(i) We compute the Hilbert polynomial from the minimal free resolution
in the usual way. For example in the case $n=5$,
\[ h(t) = \binom{t+4}{4} -5 \binom{t+2}{4} 
   + 5 \binom{t+1}{4} -\binom{t-1}{4} = 5t. \]
(ii) This is a definition. See for example \cite[I, \S7]{Ha}.
\end{Proof}

\begin{Proposition}
\label{isg1nc}
Let $n=3,4,5$ and let $\phi \in X_n$. \\
(i) If $C_\phi \subset \PP^{n-1}$ is a smooth curve of genus one 
then it is a genus one normal curve of degree $n$. \\
(ii) If $C_\phi \subset \PP^{n-1}$ is a rational curve with a single node
then it is a rational nodal curve of degree $n$.
\end{Proposition}
\begin{Proof}
By Lemma~\ref{isminres} the complex $\FF_{\bullet}(\phi)$ 
is a minimal free resolution of $R/I_\phi$.
Since $I_\phi \subset I(C_\phi)$, a comparison of Hilbert polynomials
as described in Lemma~\ref{hilbpoly} shows that $C_\phi$ has degree
$d \le n$. If $C_\phi \subset \PP^{n-1}$ spans a linear subspace
of dimension $m-1$ it follows by Riemann-Roch that
$3 \le m \le d \le n$.
We must show that $m=n$. In the case $n=3$ this is already clear.
If $n=4,5$ then $C_\phi$ is defined by quadrics. This enables us
to rule out the unwanted possibilities for $(m,d)$, with the
exception of $(m,d)=(4,4)$ in the case $n=5$. This possiblity is
excluded by the following lemma.
\end{Proof}

\begin{Lemma}
Let $C \subset \PP^3$ be either a genus one normal curve or
a rational nodal curve. Then $C$ cannot be defined by the $4 \times 4$ 
Pfaffians of a $5 \times 5$ alternating matrix of linear forms on $\PP^3$.
\end{Lemma}
\begin{Proof}
Let $\phi$ be such a matrix, with vector of submaximal Pfaffians 
$P=(p_1, \ldots, p_5)$. Let $C$ be defined by quadrics $q_1,q_2$.
By Proposition~\ref{nquads} we have $\langle p_1, \ldots, p_5 \rangle
= \langle q_1 , q_2 \rangle $. Replacing $\phi$ by $A^T \phi A$ for suitable
$A \in \GL_5$ we may suppose that $P=(q_1,q_2,0,0,0)$. 
Since $P \phi = 0$, and $q_1,q_2$ are
coprime, it follows that the first two rows of $\phi$ are zero. But then
every $4 \times 4$ Pfaffian of $\phi$ vanishes, which is a contradiction.
\end{Proof}

% \begin{Remark}
% The hypothesis of Proposition~\ref{isg1nc}(i) is that $C_\phi$ is a smooth
% curve of genus one. It would not be enough to assume that
% $C_\phi$ is a smooth curve, as the following examples with
% $C_\phi \isom \PP^1$ show
% \[ \begin{array}{lrcl} \medskip
% n=3 \qquad & \phi & = & \begin{pmatrix} x_1^3 \end{pmatrix} \\ \medskip
% n=4 & \phi & = & \begin{pmatrix} x_1^2 \\ x_2^2 \end{pmatrix} \\ 
% n=5 & \phi & = & \begin{pmatrix} 
% 0 & 0 & 0 & x_1 & x_2 \\ 
%   & 0 & x_1 & x_2 & x_3 \\
%   &   & 0 & x_3 & 0 \\
%   & -  &   & 0 & 0 \\
%   &   &   &   & 0 
% \end{pmatrix}. 
% \end{array} \]
% \end{Remark}

\begin{Lemma}
\label{israd}
Let $n=3,4,5$ and let $\phi \in X_n$. If $C_\phi \subset \PP^{n-1}$
is either a genus one normal curve or a rational nodal curve
then $I_\phi$ is a radical ideal, equivalently
$I(C_\phi) = I_\phi$. 
\end{Lemma}

\begin{Proof}
By Lemma~\ref{isminres} the complex $\FF_{\bullet}(\phi)$ 
is a minimal free resolution of $R/I_\phi$.
If $n=3$ then $I_\phi=(f)$ where $f$ is an irreducible cubic.
If $n=4,5$ then $I_\phi$ is generated by a vector space
of quadrics of dimension $d=2,5$. Since
$I_\phi \subset I(C_\phi)$ we are done by Proposition~\ref{nquads}.
\end{Proof}

\begin{Lemma}
\label{eqlem}
Let $n=3,4,5$ and let $\phi,\phi' \in X_n$. Suppose that \\
(i) there exists $\alpha \in \PGL_n$ with $\alpha(C_\phi)=C_{\phi'}$, \\
(ii) $\FF_\bullet(\phi)$ and $\FF_\bullet(\phi')$ are minimal
free resolutions of $R/I_\phi$ and $R/I_{\phi'}$, \\
(iii) the ideals $I_{\phi}$ and $I_{\phi'}$ are radical ideals. \\
Then $\phi$ and $\phi'$ are equivalent. 
\end{Lemma}

\begin{Proof}
By (i) we may assume $C_\phi = C_{\phi'}$. 
Then (iii) gives $I_{\phi} = I_{\phi'}$. 
The cases $n=3,4$ are now clear.
% If $n=3$ then $\phi' = [\mu,I_3] \phi$ for some $\mu \in K^\times$.
% If $n=4$ then $\phi' = [A,I_4] \phi$ for some $A \in \GL_2$. 
If $n=5$ then there is an isomorphism of complexes
\[ \xymatrix{ 
0 \ar[r] & R(-5) \ar[r]^-{P^T} \ar[d]^{c} 
& R(-3)^5 \ar[r]^{\phi} \ar[d]^{B} & R(-2)^5 
\ar[r]^-{P} \ar[d]^{A} & R \ar[r] \ar@{=}[d]& 0 \\ 
0 \ar[r] & R(-5) \ar[r]^-{P'^T} & R(-3)^5 \ar[r]^{\phi'} 
& R(-2)^5 \ar[r]^-{P'}   & R \ar[r] & 0 
} \]
The matrices $A,B \in \GL_5$ are uniquely determined. Comparing
this diagram with its dual gives $B = c A^{-T}$. 
So $\phi' = [A , c^{-1} I_5] \phi$. 
\end{Proof}

\begin{Lemma}
\label{stablem}
Let $n=3,4,5$ and let $\phi \in X_n$. Suppose that \\
(i) there are only finitely many $\alpha \in \PGL_n$ 
with $\alpha(C_\phi)=C_{\phi}$, \\
(ii) $\FF_\bullet(\phi)$ is a minimal free resolution of $R/I_\phi$, \\
(iii) the ideal $I_{\phi}$ is a radical ideal. \\
Then the stabiliser of $\phi$ for the action of $G_n$ 
on $\PP(X_n)$ is finite.
\end{Lemma}

\begin{Proof}
This is clear for $n=3,4$. In the case $n=5$ it suffices
to show that if $[A,I_5]\phi = \lambda \phi$ for 
some $A \in \GL_5$ and $\lambda \in K^*$, then $A$ is a scalar matrix.
Taking submaximal Pfaffians we obtain $P \adj A = \lambda^2 P$.
% For $n=4,5$ it suffices to show that if
% $[A,I_n]\phi = \lambda \phi$ for some $\lambda \in K^\times$ then $A$ is a
% scalar matrix. This is clear for $n=4$. 
% % then there are only finitely many possibilities for $A$. 
% % If $n=4$ then $A=I_2$. 
% If $n=5$ then taking submaximal Pfaffians we obtain 
% $P \adj A = \lambda^2 P$. 
By (ii) the components of $P$ are linearly
independent. It follows that $\adj A$ and hence $A$ is a scalar matrix. 
%  and hence $A = \pm I_5$. 
\end{Proof}

\begin{ProofOf}{Proposition~\ref{PropB}} (For $n=3,4,5$.)
We are given $\phi,\phi' \in X_n$ with $C_\phi$ and $C_{\phi'}$
either smooth curves of genus one or rational curves with a single node.
By Proposition~\ref{isg1nc} these are either genus one normal curves
or rational nodal curves. The hypotheses of Lemmas~\ref{eqlem} 
and~\ref{stablem} are satisfied by Lemmas~\ref{trans}, \ref{finite},
\ref{isminres} and~\ref{israd}.
\end{ProofOf}

\subsection{The generic model}
\label{sec:gen}

We show that the generic genus one model of degree $n=3,4,5$
defines a smooth curve of genus one. 

\begin{Definition}
Let $n=3,4,5$. The Jacobian matrix $J_\phi$ of a genus 
one model $\phi \in X_n$ is 
\[ \begin{array}{lll}
n=3 \quad & \phi = (f) & J_\phi = \big(\frac{\partial f}{\partial x_j}\big) \\
n=4 & \phi = \begin{pmatrix} q_1 \\ q_2 \end{pmatrix} 
& J_\phi =  \big(\frac{\partial q_i}{\partial x_j} \big) \\
n=5 & p_i = (-1)^{i+1} \pf (\phi^{\{i\}}) \quad 
& J_\phi =  \big(\frac{\partial p_i}{\partial x_j} \big). 
\end{array} \]
\end{Definition}

\begin{Lemma}
\label{myjaccrit} 
Let $n=3,4,5$ and let $\phi \in X_n$.  \\
(i) If $P \in C_\phi$ then $\rank J_\phi(P) \le n-2$. \\
(ii) If $\rank J_\phi(P) = n-2$ for every $P \in C_\phi$ then
$C_\phi$ is a smooth curve of genus one.
\end{Lemma}
\begin{Proof}
(i) We saw in Lemma~\ref{isminres}(i) that every component of $C_\phi$ has
dimension at least 1. Therefore $\dim T_P (C_\phi) \ge 1$. 
Since $I_\phi \subset I(C_\phi)$ it follows that $\rank J_\phi(P) \le n-2$. \\
(ii) The argument used in (i) shows that every component of $C_\phi$
has dimension 1. So by Lemma~\ref{isminres}(ii) 
the complex $\FF_{\bullet}(\phi)$ is a minimal free 
resolution of $R/I_\phi$. In particular $R/I_\phi$ is Cohen-Macaulay.
It follows by Serre's criterion (see \cite[\S18.3]{E}) that $I_\phi$
is a prime ideal. Hence $C_\phi$ is an irreducible smooth curve and
$I(C_\phi)= I_\phi$. It only remains to check that $C_{\phi}$ has genus $1$.
We do this by computing the Hilbert polynomial as described
in Lemma~\ref{hilbpoly}.
\end{Proof}
We define some ``bad'' subsets $B_n \subset X_n$.
\begin{Definition} (i) Let $B_3 \subset X_3$ consist of
all models of the form
  \[ \phi= \begin{pmatrix} x_1 f_1(x_2,x_3) + f_2(x_2,x_3) \end{pmatrix}. \] 
(ii) Let $B_4 \subset X_4$ consist of all models of the form 
\[ \phi = \begin{pmatrix} x_1 x_2 + g_1(x_2,x_3,x_4) \\
g_2(x_2,x_3,x_4)  \end{pmatrix}. \]
(iii) Let $B_5 \subset X_5$ consist of all models $\phi$ 
with $\phi_{ij}(1,0,0,0,0)=0$ for all $\{i,j\} \not= \{1,2\}$, 
and $\phi_{45}(x_1, \ldots, x_5) \equiv 0$. 
\end{Definition}

\begin{Lemma} 
\label{singtfae}
Let $n=3,4,5$ and let $\phi \in X_n$.  The following are equivalent. \\
(i) $C_\phi$ is not a smooth curve of genus one. \\
(ii) $\rank J_\phi(P) < n-2$ for some $P \in C_\phi$. \\
(iii) $\phi$ is equivalent to a model in $B_n$. 
% A genus one model $\phi \in X_n$ is singular if and only if 
% it is equivalent to a model in $B_n$. 
\end{Lemma}

\begin{Proof}
(i) $\Rightarrow$ (ii). This is a restatement 
of Lemma~\ref{myjaccrit}. \\
(ii) $\Rightarrow$ (i). This follows from Proposition~\ref{isg1nc}(i), 
Lemma~\ref{israd} and the Jacobian criterion for smoothness. \\
(iii)  $\Rightarrow$ (ii). Without loss of generality $\phi \in B_n$. 
Then the point $P=(1:0: \ldots :0)$ belongs to $C_\phi$ and
$\rank J_\phi(P)< n-2$. \\ %  for $P=(1:0: \ldots :0)$. \\
(ii)  $\Rightarrow$ (iii). This is clear for $n=3,4$. 
We take $n=5$.
Since $P \in C_\phi$ the $4 \times 4$ Pfaffians of $\phi(P)$ 
vanish. So $\rank \phi(P) = 0$ or $2$. 

If $\rank \phi(P) = 2$ then we may assume $P=(1:0: \ldots :0)$ and 
\[ \phi = \begin{pmatrix}
0 & x_1 & \phi_{13} & \phi_{14} & \phi_{15} \\
  &  0  & \phi_{23} & \phi_{24} & \phi_{25} \\
  &     &     0     & \ell_3 & -\ell_2 \\
  &  -  &          &  0  &  \ell_1 \\
  &     &          &     &  0 
\end{pmatrix} \]
for some $\phi_{ij}, \ell_k \in \langle x_2, x_3,x_4,x_5 \rangle$. 
Since $\rank J_\phi(P) < 3$ the linear forms $\ell_1, \ell_2, \ell_3$ 
are linearly dependent. Replacing $\phi$ by $A \phi A^T$ for suitable
$A \in \GL_5$ we may suppose that $\ell_1=0$. Then $\phi \in B_5$ as
required. 

If $\rank \phi(P)=0$ then we may assume
\[ \phi = \begin{pmatrix}
0 &  \phi_{12} & \phi_{13} & \phi_{14} & \ell_1 \\
  &  0  & \phi_{23} & \phi_{24} & \ell_2 \\
  &     &     0     & \phi_{34} & \ell_3 \\
  &  -  &          &  0  &  \ell_4 \\
  &     &          &     &  0 
\end{pmatrix} \]
for some $\phi_{ij} \in \langle x_2, x_3,x_4,x_5 \rangle$ and
$\ell_j \in   \langle x_3,x_4,x_5 \rangle$. It is clear that
$\ell_1, \ldots, \ell_4$ are linearly dependent.
Replacing $\phi$  by $A \phi A^T$ for suitable
$A \in \GL_5$ we may suppose that $\ell_4=0$. Then $\phi \in B_5$ as
required. 
\end{Proof}

\begin{ProofOf}{Proposition~\ref{PropA}} (For $n=3,4,5$.)
Let $X_n^{\sing}$ be the set of all models $\phi \in X_n$ 
which do not define a smooth curve of genus one. We consider
the projective variety 
\[ Z_n = \{ (\phi,P) \in \PP(X_n) \times \PP^{n-1} | P \in C_\phi
\text{ and } \rank J_\phi(P) < n-2 \}. \]
Let $\pr_1 : Z_n \to \PP(X_n)$ be the first projection.
Lemma~\ref{singtfae} identifies $\pr_1(Z_n)= \PP(X_n^{\sing})$.
%  where $\pr_1$ is the first projection. 
Since the image of a projective variety is again projective
it follows that $X_n^{\sing} \subset X_n$ is a Zariski closed subset. 

Lemma~\ref{singtfae} also identifies $X_n^{\sing}$ as 
the image of a morphism
\[ \G_n \times B_n \to X_n. \]
Since $\G_n$ and $B_n$ are irreducible it follows that $X_n^{\sing}$
is irreducible.
\end{ProofOf}

\subsection{The invariant differential}
\label{sec:invdiff}
We continue to work over an algebraically closed field $K$.
In the case $n=5$ we further suppose that $\charic(K)  \not=2$.
% In this subsection we further suppose that $\charic(K) \not=2,3$.

Let $\phi \in X_n$ with $C_\phi$ a smooth of curve genus one.
We use $\phi$ to define an invariant differential
$\omega_\phi$ on $C_\phi$. 
% We show that $\phi$ determines a non-zero regular
% 1-form on $\omega_\phi$ on $C_\phi$. 
In the cases $n=1,2$ we put
\[ \begin{array}{lll} \smallskip
n = 1 \quad & \phi = (a_1,a_2,a_3,a_4,a_6) \quad & 
\omega_\phi   =  \frac{ dx }{ 2 y + a_1 x + a_3  } \\
n = 2 & \phi = (p(x,z),q(x,z)) &  
\omega_\phi   =  \frac{z^2 d(x/z)}{2y+p(x,z)}.
\end{array} \]
In the cases $n=3,4,5$ we start with the complex 
\[ \FF_{\bullet}(\phi) \, : \quad
0 \ra R \stackrel{\phi_{n-2}}{\ra} \FF_{n-3}
\ra \ldots \ra \FF_1 \stackrel{\phi_1}{\ra} R \ra 0 \]
defined in \S\ref{minfreeres}. We identify the maps $\phi_i$ 
with the matrices of homogeneous polynomials that represent them. 
Then we define
\[ \omega_\phi = \frac{ x_1^2 d(x_2/x_1)}
{ \frac{\partial \phi_1}{\partial x_3}
\circ \ldots \circ \frac{\partial \phi_{n-2}}{\partial x_{n}} } \]
where the partial derivative of a matrix is 
the matrix of partial derivatives.
In the cases $n=3,4$ this formula works out as
\[ \omega_\phi = \frac{x_1^2 d(x_2/x_1)}{\frac{\partial f}{\partial x_3}} 
\quad \text{ and } \quad  \omega_\phi = 
\frac{x_1^2 d(x_2/x_1)}{
\frac{\partial q_1}{\partial x_4} \frac{\partial q_2}{\partial x_3} -
\frac{\partial q_1}{\partial x_3} \frac{\partial q_2}{\partial x_4}}. \]

\begin{Proposition}
\label{omegaprop}
Let $\phi \in X_n$ with $C_\phi$ a smooth curve of genus one. 
% Let $g \in \G_n$ and put $\phi' = g \phi$.
If $\phi' = g \phi$ for some $g \in \G_n$ then the isomorphism 
$\gamma : C_{\phi'} \isom C_\phi$ determined by $g$ satisfies  
\[ \gamma^* \omega_\phi = (\det g) \omega_{\phi'}. \]
\end{Proposition}
\begin{Proof}
If the proposition holds for $g_1,g_2 \in \G_n$ then it holds for $g_1 g_2$. 
So we only need to consider $g$ running over
a set of generators for $\G_n$. Since the cases $n=1,2$ are well known
we take $n=3,4,5$. The result is clear for $g$ of the form $[1,B]$ with
\[ B = \begin{pmatrix} 
 \mu_1 & \lambda \\
& \mu_2 \\
& & \ddots \\
& & & \mu_n
\end{pmatrix}. \]

If $n=3$ and $g=[\mu,I_3]$ then the result is again clear.
If $n=4$ and $g=[A,I_4]$ then there is an isomorphism of complexes
\[ \xymatrix{ 
0 \ar[r] & R(-4) \ar[r]^{\phi'_2} \ar[d]^{\det A} & R(-2)^2 
\ar[r]^-{\phi'_1} \ar[d]^{A^T} & R \ar[r] \ar@{=}[d]& 0 \\ 
0 \ar[r] & R(-4) \ar[r]^{\phi_2} 
& R(-2)^2 \ar[r]^-{\phi_1}   & R \ar[r] & 0 
} \]
We deduce 
\[\gamma^* \omega_\phi = (\det A) \omega_{\phi'}
= (\det g) \omega_{\phi'}. \]

If $n=5$ and $g=[A,I_5]$ then there is an isomorphism of complexes
\[ \xymatrix{ 
0 \ar[r] & R(-5) \ar[r]^-{\phi'_3} \ar[d]^{(\det A)^2} 
& R(-3)^5 \ar[r]^{\phi'_2} \ar[d]^{(\det A) A^T} & R(-2)^5 
\ar[r]^-{\phi'_1} \ar[d]^{\adj A} & R \ar[r] \ar@{=}[d]& 0 \\ 
0 \ar[r] & R(-5) \ar[r]^-{\phi_3} & R(-3)^5 \ar[r]^{\phi_2} 
& R(-2)^5 \ar[r]^-{\phi_1}   & R \ar[r] & 0 
} \]
We deduce 
\[\gamma^* \omega_\phi = (\det A)^2 \omega_{\phi'}
= (\det g) \omega_{\phi'}. \]

It only remains to prove the proposition for $g=[1,B]$ with $B$ 
a permutation matrix. This in turn reduces to checking the 
result for a set of transpositions generating the 
symmetric group $S_n$. The symmetry (12) is already clear
from the identity
\[ x_1^2 d (x_2/x_1) + x_2^2 d ( x_1 / x_2) = 0. \]
Since the entries of $\phi_1$ belong to $I(C_{\phi})$ we have
\begin{equation}
\label{omega*} 
\sum_{i=2}^{n} \frac{\partial \phi_1}{\partial x_i} d(x_i/x_1)  =  0. 
\end{equation}
If $n=3$ then  (\ref{omega*}) gives the symmetry  $(23)$. 
If $n=4$ then the symmetry $(34)$ is clear. By~(\ref{omega*}) we have
\[ \left| \begin{matrix} \smallskip
\frac{\partial q_1}{\partial x_2} & \frac{\partial q_1}{\partial x_4} \\
\frac{\partial q_2}{\partial x_2} & \frac{\partial q_2}{\partial x_4} 
\end{matrix} \right|  d( x_2/x_1)  +
\left| \begin{matrix} \smallskip
\frac{\partial q_1}{\partial x_3} & \frac{\partial q_1}{\partial x_4} \\
\frac{\partial q_2}{\partial x_3} & \frac{\partial q_2}{\partial x_4} 
\end{matrix} \right|  d( x_3/x_1) =0  \]
and this establishes the symmetry $(23)$. 

If $n=5$ then the symmetry $(35)$ is clear.
Differentiating $\phi_1 \phi_2 = 0$ and $\phi_2 \phi_3= 0$ we find
\begin{equation}
\label{diffid}
 \frac{\partial \phi_1}{\partial x_3} 
   \frac{\partial \phi_2}{\partial x_4} 
   \frac{\partial \phi_3}{\partial x_5} 
 +
   \frac{\partial \phi_1}{\partial x_3} 
   \frac{\partial \phi_2}{\partial x_5} 
   \frac{\partial \phi_3}{\partial x_4} 
= 
   \phi_1 
   \frac{\partial \phi_2}{\partial x_3} 
   \frac{\partial^2 \phi_3}{\partial x_4 \partial x_5}. 
\end{equation}
% Since the entries of $\phi_1$ vanish on the curve, 
This establishes the symmetry $(45)$.
Using~(\ref{omega*}) we get
\[ \sum_{i=2}^5  \frac{\partial \phi_1}{\partial x_i} 
   \frac{\partial \phi_2}{\partial x_4} 
   \frac{\partial \phi_3}{\partial x_5} 
 d( x_i/x_1 ) = 0. \]
The terms for $i=4,5$ vanish since $\charic(K) \not=2$
and analogous to (\ref{diffid}) we have 
% $\phi_2$ is alternating,
% $\phi_3 = \phi_1^T$, and analogous to (\ref{diffid}) we have 
%  The term for $i=4$ vanishes 
% Notice that the terms for $i=1,2$ vanish, the latter
% since $\charic(K) \not= 2$ and analogous to (\ref{diffid}) we have 
\begin{equation*}
  2 \, \frac{\partial \phi_1}{\partial x_4} 
   \frac{\partial \phi_2}{\partial x_4} 
   \frac{\partial \phi_3}{\partial x_5} 
= 
  \frac{\partial^2 \phi_1}{\partial x_4^2 } 
    \frac{\partial \phi_2}{\partial x_5} 
  \phi_3  \, . 
\end{equation*}
% This is where we use our assumption $\charic(K) \not= 2$.
This establishes the symmetry $(23)$.
\end{Proof}

\begin{Lemma}
\label{weieromega}
Let $\phi = \pi_n(\phi_1)$ be a Weierstrass model with $C_\phi$ 
a smooth curve of genus one. Then the natural isomorphism 
$\gamma : C_{\phi_1} \isom C_{\phi} $ satisfies
$\gamma^* \omega_{\phi} = \omega_{\phi_1}. $ 
\end{Lemma}
\begin{Proof}
We check this by direct calculation using the definition of 
$\omega_\phi$ and the formulae of \S\ref{sec:wmodels}.
\end{Proof}

\begin{Lemma}
Let $\phi \in X_n$ with $C_\phi$ a smooth curve of genus one. 
Then $\omega_\phi$ is an invariant differential on $C_\phi$.
\end{Lemma}
\begin{Proof}
Our claim is that $\omega_\phi$ is a non-zero regular 1-form.
By Propositions~\ref{PropB}, \ref{PropC} and~\ref{omegaprop} it suffices to 
prove this for $\phi$ a Weierstrass model.
Then Lemma~\ref{weieromega} reduces us to the case $n=1$,
and in this case the result is well known. 
\end{Proof}

We recall from Definition~\ref{def:propeq} 
that models $\phi, \phi' \in X_n$ are properly equivalent if 
there exists $g \in \G_n$ with $g \phi = \phi'$ and $\det(g)=1$. 

\begin{Lemma}
\label{reduceto1}
Let $\phi,\phi' \in X_1$ with $C_\phi$ and $C_{\phi'}$ smooth
curves of genus one. Then $\phi$ and $\phi'$ are properly
equivalent if and only if $\pi_n(\phi)$ and $\pi_n(\phi')$ are 
properly equivalent.
\end{Lemma}
\begin{Proof}
One implication is clear from Proposition~\ref{PropC}. For 
the converse we suppose $\pi_n(\phi)$ and $\pi_n(\phi')$ are 
properly equivalent. Then by Proposition~\ref{omegaprop} 
and Lemma~\ref{weieromega} there is an isomorphism 
$\gamma : C_\phi \isom C_{\phi'}$ with
$\gamma^* \omega_{\phi'} = \omega_\phi$. 
Composing with a translation % (if necessary)
we may suppose that $\gamma$ is determined by some
$g \in \G_1$. It follows that $\phi$ and $\phi'$ are properly equivalent.
% as was to be shown.
\end{Proof}

\medskip

\begin{ProofOf}{Proposition~\ref{PropD}}
Let $\phi \in X_n$ with $C_\phi$ a smooth curve of genus one. 
We must show that $\phi$ is properly equivalent to a Weierstrass model
$\pi_n(0,0,0,A,B)$ for some unique $A,B \in K$. 
The existence is already clear from Propositions~\ref{PropB}
and~\ref{PropC}. To prove uniqueness we use Lemma~\ref{reduceto1} to
reduce to the case $n=1$. In this case the result is well known.
\end{ProofOf}

In the proof of Theorem~\ref{mainthm}(iii) we used
\begin{Proposition}
\label{invginv}
Let $\phi \in X_n$ with $C_\phi$ a smooth curve of genus one. 
Then the geometric invariants of $(C_\phi,\omega_\phi)$ 
are $c_4(\phi)$ and $c_6(\phi)$.
\end{Proposition}
\begin{Proof}
We are free to replace $\phi$ by any properly equivalent model.
So we may assume that $\phi$ is a Weierstrass model. Then
Lemma~\ref{weieromega} reduces us to the case $n=1$. 
In this case the result is a tautology.
\end{Proof}

\section{Weierstrass models}
\label{sec:wmodels}

Let $E$ be an elliptic curve with Weierstrass equation
\[ y^2 + a_1 x y + a_3 y = x^3 + a_2 x^2 + a_4 x + a_6. \]
In the notation of \S\ref{sec:g1models} we have $E=C_\phi$ where 
$\phi = (a_1,a_2,a_3,a_4,a_6)$. %  \in X_1$.
The complete linear system $|n.0|$ determines a morphism
$E \to \PP^{n-1}$ 
\[ \begin{array}{llcl}
n=2 \qquad  & (x,y) & \mapsto & (x:1) \\
n=3 & (x,y) & \mapsto & (1:x:y)  \\
n=4 & (x,y) & \mapsto & (1:x:y:x^2) \\
n=5 & (x,y) & \mapsto & (1:x:y:x^2:xy).
\end{array} \]
The image is defined by a genus one model
\[ \begin{array}{rcl} \medskip
\pi_2(\phi) & = & 
% \begin{pmatrix}
(a_1 x z + a_3 z^2 , 
x^3 z + a_2 x^2 z^2 + a_4 x z^3 + a_6 z^4)
% \end{pmatrix} 
\\ \medskip
\pi_3(\phi) & = &
% \begin{pmatrix}
(y^2 z + a_1 x y z + a_3 y z^2 - x^3 - a_2 x^2 z - a_4 x z^2 - a_6 z^3 ) 
% \end{pmatrix}
\\ \medskip
\pi_4(\phi) & = & \begin{pmatrix} 
x_1 x_4 - x_2^2 \\
x_3^2  + a_1 x_2 x_3 + a_3 x_1 x_3 - x_2 x_4 - a_2 x_2^2 
- a_4 x_1 x_2 - a_6 x_1^2
\end{pmatrix} \\  
\pi_5(\phi) & = &  \begin{pmatrix}
0 & \ell & x_5 & x_4 & x_3 \\
 & 0 & x_4 & x_3 & x_2 \\
 & & 0 & -x_2 & 0 \\
 & - & & 0 & x_1 \\
 & & & & 0
\end{pmatrix}
\end{array} \]
where $\ell = a_1 x_5  - a_2 x_4 + a_3 x_3 - a_4 x_2  -a_6 x_1$.

These formulae define a morphism $\pi_n : X_1 \to X_n$. 
% The models
% $\pi_n(\phi)$ for $\phi \in X_1$ are called Weierstrass models.
% We refer to them collectively as the Weierstrass family.
A morphism $\gamma_n : \G_1 \to \G_n$ with the
properties specified in Proposition~\ref{PropC} is given by
\[ \begin{array}{rcl} \medskip
\gamma_2([u;r,s,t]) & = & \left[u^{-3}, (0,u^2 s,t) , 
\begin{pmatrix}  u^2 & 0 \\ r & 1 \end{pmatrix} \right] \\ \medskip
\gamma_3([u;r,s,t]) & = & \left[u^{-6}, \begin{pmatrix} 
 1 & r & t \\ 0 & u^2 & u^2 s \\ 0 & 0 & u^3  \end{pmatrix} \right] \\
% u^2 & u^2 s & 0 \\ 0 & u^3 & 0  \\ r & t & 1  \end{pmatrix} \right] \\
\gamma_4([u;r,s,t]) & = & \left[ \begin{pmatrix}
u^{-4} & 0 \\ u^{-6} r & u^{-6} \end{pmatrix} ,
\begin{pmatrix} 
1 & r & t & r^2  \\
0 & u^2 & u^2 s & 2u^2 r  \\ 
0 & 0 & u^3 & 0 \\  
0 & 0 & 0 & u^4
\end{pmatrix} \right] 
\end{array} \]
and $\gamma_5([u;r,s,t]) = \left[ A_5,B_5 \right]$
where
\[ A_5 = u^{-2} \begin{pmatrix}
1 & - s  &  2 r - s^2 & r s - t & -r^2 + r s^2 - s t \\
0 & u & 2 u s & -u r & u (-2 r s + t) \\
0 & 0 & u^2 & 0 & -u^2 r \\
0 & 0 & 0 & u^3 & u^3 s \\
0 & 0 & 0 & 0 & u^4 
 \end{pmatrix} \]
and
\[ B_5 = 
u^{-3} \begin{pmatrix} 
1 & r & t & r^2 & r t \\
0 & u^2 & u^2 s & 2 u^2 r & u^2 (r s + t) \\
0 & 0 & u^3 & 0 & u^3 r \\
0 & 0 & 0 & u^4 & u^4 s \\
0 & 0 & 0 & 0 & u^5
\end{pmatrix}.  \]
% In fact we had to choose $\pi_5$ rather
% carefully to ensure that $\gamma_5$ does not 
% depend on the $a_i$. 
% It is routinely checked that $\gamma_n$ has the required properties.
% group homomorphism making the required diagram commute. 
% Finally, in each case we have
% \[\det(\gamma_n([u;r,s,t])) = u^{-1} = \det([u;r;s,t]).\]

\section{Formulae}
\label{sec:formulae}

We recall some formulae for the invariants in the cases 
$n=2,3,4$. In each case we scale the invariants so as to give
the usual formulae when restricted to the Weierstrass family.
As noted in Remark~\ref{primrem} these are also the scalings, unique up
to sign, for which $c_4$, $c_6$ and $\Delta$ are primitive integer
coefficient polynomials. 
% Our scalings are in general different 
% from those found in the classical literature. 
% The author has contributed these formulae to {\sf MAGMA} \cite{magma}.
We assume for simplicity that $\charic(K) \not=2,3$.

% The formulae have 
% been implemented in {\sf MAGMA} \cite{magma} by the author.
% We scale the invariants to agree with the usual formulae when
% restricted to the Weierstrass family.
% Classically one includes multinomial coefficients
% in the expressions for $f$ and $U$. This has effect of making
% the coefficients in the expressions for $c_4$ and $c_6$ smaller.
% Our conventions are the right ones for number theory.
% It has all been implemented in {\sf MAGMA} by the author.

\subsection{Formulae for the invariants: case $n=2$}
The invariants in the case $n=2$ are classical. Here is
one way to compute them. We start with the binary quartic 
\[f  = a x^4 + b x^3 z + c x^2 z^2 + d x z^3 + e z^4 \]
and compute (a scalar multiple of) its Hessian
\[ \begin{array}{rcl} H & = & 
% \frac{1}{3} \det (\frac{\partial^2 f}{\partial x_i,\partial x_j}) \\ & = & 
(8 a c - 3 b^2) x^4 + (24 a d - 4 b c) x^3 z 
+ (48 a e + 6 b d - 4 c^2) x^2 z^2  \\ &  &
\, \quad \,\,  + (24 b e - 4 c d) x z^3 + (8 c e - 3 d^2) z^4. 
\end{array} \]
We then turn $f$ into a differential operator by substituting 
$\partial/\partial z$ and $-\partial/\partial x$ for $x$ and $z$.
Letting this operator act on $f$ and $H$ gives the invariants
\[ \begin{array}{rcl}
c_4 & = & 2^4(12 a e - 3 b d + c^2) \\
c_6 & = & 2^5(72 a c e - 27 a d^2 - 27 b^2 e + 9 b c d - 2 c^3).
\end{array} \]
The discriminant $\Delta = (c_4^3-c_6^2)/1728$ is $16$ times the 
usual discriminant of a degree $4$ polynomial.
If the cross terms $\alpha_0, \alpha_1, \alpha_2$ are included
(by computing the square) then
$c_4$, $c_6$ and $\Delta$ are primitive integer coefficient polynomials
in $\alpha_0, \alpha_1, \alpha_2,a,b,c,d,e$.

\subsection{Formulae for the invariants: case $n=3$}
The invariants in the case $n=3$ are again classical. 
The ternary cubic
\[ \begin{array}{rcl} 
U(x,y,z) & = & a x^3 + b y^3 + c z^3 + a_2 x^2 y + a_3 x^2 z \\
& & ~\qquad \qquad ~+ b_1 x y^2 + b_3 y^2 z + c_1 x z^2  + c_2 y z^2 + m x y z 
\end{array} \]
has Hessian 
\[ H(U) = (-1/2) \times  \left| \begin{array}{ccc} \smallskip
\frac{\partial^2 U}{\partial x^2} &
\frac{\partial^2 U}{\partial x \partial y} &
\frac{\partial^2 U}{\partial x \partial z} \\ \smallskip
\frac{\partial^2 U}{\partial x \partial y} &
\frac{\partial^2 U}{\partial y^2} &
\frac{\partial^2 U}{\partial y \partial z} \\
\frac{\partial^2 U}{\partial x \partial z} &
\frac{\partial^2 U}{\partial y \partial z} &
\frac{\partial^2 U}{\partial z^2} 
\end{array} \right| \enspace. \]
% \[ H(U) = (-1/2) \det 
% \bigg( \frac{\partial^2 U}{\partial x_i \partial x_j} \bigg).\] 
% The factor $-1/2$, although not standard, is a choice we find convenient.
% The Hessian is a polynomial map $H : X_3 \to X_3$ satisfying
% \[ H \circ g = (\det g)^2 g \circ H\] 
% for all $g \in \GL_3(\Kbar)$. We say it is 
% covariant of weight $2$. 
Putting $c_4= c_4(U)$, $c_6= c_6(U)$ and $H=H(U)$ we find
\[H( \lambda U + \mu H) =3 (c_4 \lambda^2 \mu + 2 c_6 \lambda \mu^2 
 + c_4^2 \mu^3) U 
 + (\lambda^3 - 3 c_4  \lambda \mu^2- 2 c_6 \mu^3  ) H. \]
This formula is classical: see \cite[\S II.7]{Hilbert} or 
\cite[\S225]{Salmon}. It is easily
verified by restricting to any family of plane cubics covering
the $j$-line, for example the Weierstrass family 
defined in \S\ref{sec:wmodels}.
We solve to find
\[ \begin{array}{rcl} c_4 & = &
-216 a b c m + 144 a b c_1 c_2 + 144 a c b_1 b_3 - 48 a b_1 c_2^2 
+ \ldots \\ & & 
~\quad \ldots + 24 a_3 b_1 c_2 m - 8 a_3 b_3 m^2 + 16 b_1^2 c_1^2 - 8 b_1 c_1 m^2 + m^4 \medskip \\ 
c_6 & = &
5832 a^2 b^2 c^2 - 3888 a^2 b c b_3 c_2 + 864 a^2 b c_2^3 + 864 a^2 c b_3^3 +
\ldots \\ & & 
~\quad \ldots  + 12 a_3 b_3 m^4 + 64 b_1^3 c_1^3 - 
48 b_1^2 c_1^2 m^2 + 12 b_1 c_1 m^4 - m^6
\end{array}
\]
where the full expressions have 25 terms and 103 terms respectively.
These polynomials are written out completely in 
\cite{Mc+}, \cite[\S10.3]{Dolg}, 
\cite[\S\S220,221]{Salmon} and \cite[\S\S 4.4,4.5]{St}.

We may compute the discriminant as $\Delta = (c_4^3-c_6^2)/1728$. 
An alternative, taken from \cite[\S241]{Salmon}, is the following.
We compute the partial derivatives of $U$ and $H=H(U)$, and 
arrange the coefficients of these quadrics in a $6 \times 6$ matrix.
Then this matrix has determinant $\pm 1728 \Delta$.

\subsection{Formulae for the invariants: case $n=4$}

We identify a genus one model of degree $4$ with a pair of 
$4 \times 4$ symmetric matrices. Explicitly
\[ \phi = \begin{pmatrix} q_1 \\ q_2 \end{pmatrix} 
\equiv \begin{pmatrix} A \\ B \end{pmatrix} \]
where
% \[ q_1(x_1, \ldots x_4) = \frac{1}{2} \x^T A \x \quad \text{ and } \quad
% q_2(x_1, \ldots x_4) = \frac{1}{2} \x^T B \x. \]
\[ \begin{array}{rcl} q_1(x_1, \ldots, x_4) & = & \frac{1}{2} \x^T A \x 
\end{array} \,\, \text{ and } \,\, \begin{array}{rcl}
q_2(x_1, \ldots, x_4) & = & \frac{1}{2} \x^T B \x. \end{array} \]
The invariants are found by computing the binary quartic
\[ \det(sA+tB) = a s^4 + b s^3 t + c s^2 t^2 + d s t^3+e t^4 \]
and then using the formulae for $n=2$. The correct scalings are
\[ \begin{array}{rcl}
c_4 & = & 12 a e - 3 b d + c^2 \\
c_6 & = & \frac{1}{2} (72 a c e - 27 a d^2 - 27 b^2 e + 9 b c d - 2 c^3). 
\end{array} \]
Since $b,d \in 2 \Z[X_4]$ the coefficients of $c_4$ and $c_6$ are 
indeed integers as predicted by Lemma~\ref{intcoeff}

% so contrary to inital impressions
% $c_6$ does indeed have integer coefficients.
% In fact, as prediced by Proposition~??, $c_4$ and $c_6$ 
% are primitive polynomials in $\Z[X_4]$. They have $4003$ and $81547$
% terms respectively.

We may compute the discriminant as $\Delta = (c_4^3-c_6^2)/1728$. 
An alternative is the following. Let $T_1$ and $T_2$ be the symmetric
matrices defined in~\cite{Mc+},~\cite{MSS} by
\[ \begin{array}{rcl}
\adj(s(\adj A)+t(\adj B)) & = 
& a^2 A s^3 + a T_1 s^2 t + e T_2 s t^2 + e^2 B t^3.
\end{array}
\]
The corresponding quadrics are
% \[ q'_1(x_1, \ldots x_4) = \frac{1}{2} \x^T T_1 \x \quad \text{ and } \quad
% q'_2(x_1, \ldots x_4) = \frac{1}{2} \x^T T_2 \x. \]
\[ \begin{array}{rcl} q'_1(x_1, \ldots x_4) & = & \frac{1}{2} \x^T T_1 \x 
\end{array} \,\, \text{ and } \,\, \begin{array}{rcl}
q'_2(x_1, \ldots x_4) & = & \frac{1}{2} \x^T T_2 \x. \end{array} \]
For a permutation $\pi \in S_4$ we define
\[ \begin{array}{rcl} \Omega_{\pi(1),\pi(2)}  & = & \sign(\pi) \left( 
 \frac{\partial q_1}{\partial x_{\pi(3)}} 
  \frac{\partial q_2}{\partial x_{\pi(4)}}
- \frac{\partial q_1}{\partial x_{\pi(4)}} 
  \frac{\partial q_2}{\partial x_{\pi(3)}} \right). \end{array} \]
Then we arrange the coefficients of the 
quadrics $q_1$, $q_2$, $q'_1$, $q'_2$ and $ \Omega_{r,s} $ 
% = \frac{\partial q_1}{\partial x_r} \frac{\partial q_2}{\partial x_s}
% - \frac{\partial q_1}{\partial x_s} \frac{\partial q_2}{\partial x_r}\[
for $1 \le r < s \le 4$ in a $10 \times 10$ matrix. 
The determinant of this matrix turns out to be $\pm 16 \Delta$.
As seen in \S\ref{sec:invdiff}, the quadrics $\Omega_{r,s}$ arise naturally in 
the construction of an invariant differential $\omega_\phi$ on $C_\phi$.

\section{An evaluation algorithm}
\label{sec:evalalg}
% \subsection{Evaluation algorithms for $n=5$}
In the case $n=5$ the invariants $c_4$ and $c_6$ are homogeneous polynomials
of degrees $20$ and $30$ in $50$ variables. They are therefore 
too large to write down as explicit polynomials. Nonetheless 
we have found a practical algorithm for evaluating them.
%  and hence by Theorem~\ref{mainthm} a practical algorithm for computing
% the Jacobian of a genus one normal quintic.
We assume throughout this section that $\charic(K) \not= 2,3,5$.

We identify $X_5 = \wedge^2 V \otimes W$ where $V$ and $W$ are 
$5$-dimensional vector spaces. Explicitly %  we identify
\[ ( \phi_{ij}(x_1, \ldots, x_5))_{i,j = 1, \ldots,5} \equiv
\sum_{i<j}  (v_i \wedge v_j) \otimes \phi_{ij}(x_1, \ldots, x_5) \]
where $v_1, \ldots, v_5$ and $x_1, \ldots,x_5$ are fixed bases
for $V$ and $W$. The action of $\G_5 = \GL(V) \times \GL(W)$ is
the natural one. The commutator subgroup of $\G_5$ 
is $G_5= \SL(V) \times \SL(W)$.

\begin{Definition}
Let $(\rho,Y)$ be a rational representation of $\G_5$. A covariant
is a polynomial map $F : \wedge^2 V \otimes W \to Y$ such that
$F \circ g= \rho(g) \circ F$ for all $g \in G_5.$
\end{Definition}

Notice that the invariants are the covariants in the case 
of the trivial representation. For a fixed representation $(\rho,Y)$ 
the covariants form a module over the ring of invariants.

The $4 \times 4$ Pfaffians of $\phi$ are quadrics $p_1, \ldots, p_5$ 
satisfying 
\[ \phi \wedge \phi \wedge v_i = p_i(x_1, \ldots, x_5) \,\, v_1 \wedge \ldots
\wedge v_5. \] 
We may therefore define covariants
\[ \begin{array}{ll}
P : \wedge^2 V \otimes W \to V^* \otimes S^2W \, ; & \phi \mapsto 
\sum_{i=1}^5 v_i^* \otimes p_i(x_1, \ldots, x_5) \\
S: \wedge^2V \otimes W \to S^5 W \, ;  & 
  \phi \mapsto \det (\frac{\partial p_i}{\partial x_j})
\end{array} \]
where $v_1^*, \ldots, v_5^*$ is the basis for $V^*$ dual 
to $v_1, \ldots,v_5$. 

Our method for evaluating the invariants relies on the following
geometric ``accident''.
\begin{Lemma} 
\label{auxexists}
Let $\phi \in X_5$ with $C_\phi$ a smooth curve of genus one
and let $p_1 ,\ldots, p_5$ be the $4 \times 4$ Pfaffians of $\phi$. \\
% and let $s= S(\phi)$. Then \\
(i) The secant variety of $C_\phi$ is the hypersurface 
defined by $S(\phi)=0$. \\
(ii) The partial derivatives $\frac{\partial}{\partial x_i} S(\phi)$ 
are quadrics in $p_1, \ldots,p_5$. 
\end{Lemma}
\begin{Proof} 
(i) See~\cite[Lemma 6.7]{g1pf} or~\cite[VIII.2.5]{Hu}. \\
(ii) The lemma may be checked by direct computation on any family 
of models covering the $j$-line, for example the 
Weierstrass family defined in \S\ref{sec:wmodels}.
A more illuminating proof is given in~\cite[Corollary 7.5]{hsecenc}. 
\end{Proof}
Lemma~\ref{auxexists}(ii) is accompanied by the following uniqueness
statement. 
\begin{Lemma}
\label{auxunique}
Let $\phi \in X_5$ with $4 \times 4$ Pfaffians $p_1 ,\ldots, p_5$.
If $S(\phi) \not= 0$ then the quartics $\{ p_i p_j : 1 \le i \le j \le 5 \}$ 
are linearly independent.
\end{Lemma}

\begin{Proof} The condition $S(\phi) \not=0$ gives that
$p_1 ,\ldots, p_5$ are linearly independent.
Now suppose $q(v_1, \ldots, v_5)$ is a quadric in 5 variables
with $q(p_1, \ldots, p_5)= 0.$  We differentiate with 
respect to $x_j$ to obtain \[ \begin{array}{c} \sum_{i=1}^5 
\frac{\partial q}{\partial v_i}(p_1, \ldots,p_5) 
\frac{\partial p_i}{\partial x_j}(x_1, \ldots,x_5) = 0. \end{array} \]
Our assumption $S(\phi) \not=0$ then gives
$\frac{\partial q}{\partial v_i}(p_1, \ldots,p_5) = 0$
for all $i$. Since $p_1, \ldots, p_5$ are linearly independent,
it follows that all partial derivatives of $q$ 
% $\frac{\partial q}{\partial v_i}$ 
are identically zero, and hence that $q$ itself is identically zero.
\end{Proof}

\begin{Lemma} 
\label{Qcov}
There is a covariant 
\[ \begin{array}{ll}
Q : \wedge^2 V \otimes W \to S^2 V \otimes W \, ; & \phi \mapsto 
\sum q_i(v_1, \ldots, v_5) \otimes x_i 
\end{array} \]
with the property that if $\phi \in X_5$ with $4 \times 4$ 
Pfaffians $p_1, \ldots, p_5$ then 
% \[ \begin{array}{rcl}
% P(\phi) & = & \sum v_i^* \otimes p_i(x_1, \ldots, x_5) \\
% Q(\phi) & = & \sum q_i(v_1, \ldots, v_5) \otimes x_i \\
% S(\phi) & = & s(x_1, \ldots, x_5) 
% \end{array} \]
% then 
\[\begin{array}{rcl} 
\frac{\partial }{\partial x_i} S(\phi) 
% \frac{\partial S(\phi) }{\partial x_i} 
&  = & q_i(p_1, \ldots, p_5)
\end{array} \] for all $i$. 
Moreover $Q$ is uniquely determined by this property.
\end{Lemma}
\begin{Proof}
Let $\phi \in X_5(\K)$ be the generic model defined over the function 
field $\K = K(X_5)$.  By Proposition~\ref{PropA} we know that $C_\phi$ 
is a smooth curve of genus one. Then by Lemma~\ref{auxexists} we can
solve for quadrics $q_1, \ldots ,q_5$ with the required property.
These quadrics define a rational map
\[Q : \wedge^2 V \otimes W \ratto S^2 V \otimes W.\]
By Lemma~\ref{auxunique} the quadrics $q_1, \ldots, q_5$
are uniquely determined. So the covariance property is clear.
We must show that $Q$ is regular, and for this we may work over
an algebraically closed field. 

We first claim that $Q$ is regular at all $\phi \in X_5$
with $S(\phi) \not=0$. The coefficients of the quartics
$\{ p_i p_j : 1 \le i \le j \le 5 \}$ may be arranged in a $15 \times 70$
matrix. Let $h_1, \ldots, h_N \in K[X_5]$ be 
the $15 \times 15$ minors of this matrix. If $\phi \in X_5$ with 
$S(\phi) \not=0$ then Lemma~\ref{auxunique} gives $h_i(\phi) \not= 0$
for some $i$. Our claim follows since $Q$ is regular on 
each of the open sets $\{ h_i \not= 0 \}$.

Now let $F \in K[X_5]$ be a homogeneous polynomial of 
least degree such that $FQ$ is regular. 
Then $Q$ is regular at $\phi$ if and only if $F(\phi) \not=0$.
The above claim gives $F(\phi) \not= 0$ whenever $S(\phi) \not=0$.
But we know by Lemma~\ref{auxexists}(i) that 
$S(\phi) \not=0$ for $C_\phi$ a smooth curve of genus one.
By Theorem~\ref{mainthm}(ii) and the irreducibility of $\Delta$ 
(which is inherited from the case $n=1$) it
follows that $F$ is a power of $\Delta$. To complete the proof it
only remains to show that $S$ is not divisible by $\Delta$. 
Since $S$ has degree 10 and $\Delta$ has degree 60, this is clear.
\end{Proof}

% In Lemma~\ref{Qcov} we show that $Q$ extends to a covariant.
Starting from $P$ and $Q$ we compute covariants $M$ and $N_\la$ 
taking values in $S^5 V^*$ and $S^5 V$. 
We then use the natural identification $S^5 V^* = (S^5V)^*$ 
to contract these covariants, and hence compute the invariants.
We arrive at the following algorithm.

\begin{Algorithm} 
\label{myalg}
Assume $\charic(K) \not= 2,3,5.$ \\
INPUT: A genus one model $\phi \in X_5 = \wedge^2 V \otimes W$. \\
OUTPUT: The invariants $c_4(\phi)$, $c_6(\phi)$, $\Delta(\phi)$. 
\begin{enumerate}
\item Compute the $4 \times 4$ Pfaffians $p_1, \ldots, p_5$ of $\phi$. 
\item Check that the quartics $\{ p_i p_j : 1 \le i \le j \le 5 \}$ 
are linearly independent. If not return $0,0,0$. 
\item Compute the secant quintic $s=\det( \frac{\partial p_i}{\partial x_j})$.
\item Solve for the auxiliary quadrics $q_1, \ldots,q_5$ satisfying
\[\begin{array}{rcl}
\frac{\partial s}{\partial x_i} & = & q_i(p_1, \ldots ,p_5).
\end{array} \]
\item Compute the quintic $M = \det(\sum_{k=1}^5 
\frac{\partial^2 p_k}{\partial x_i \partial x_j} v_k^* ) \in S^5 V^*$. 
\item Compute the quintic $N_\lambda = 
\det ( \lambda \frac{\partial q_i}{\partial v_j} 
  + \sum_{k=1}^5 \frac{\partial \phi_{jk}}{\partial x_i} v_k  ) \in S^5 V$.
\item Contract $M$ and $N_\lambda$ to obtain
% $\parital/\partial v_i$ for $v_i^*$. Then acting on $N_t$ we obtain
\[ \langle M,N_\lambda \rangle = 40 c_4 \lambda 
  - 320 c_6 \lambda^3 + 128 c_8 \lambda^5. \]
\item Check that $c_8= c_4^2$. 
\item Return $c_4$, $c_6$, $(c_4^3-c_6^2)/1728$. 
\end{enumerate}
\end{Algorithm}
 
It is easy to show that the quantities $c_4$ and $c_6$
computed % in Algorithm~\ref{myalg} 
are invariants of weights $4$ and $6$.
By Theorem~\ref{mainthm} the invariants of
weights $4$ and $6$ each form a 1-dimensional vector space. 
So it only remains to check that the invariants 
computed are not identically zero, 
and that they are correctly scaled. We did this by computing 
their restriction to the Weierstrass family, but in fact it would
suffice to compute a single numerical example.

To complete the justification of Algorithm~\ref{myalg} we must show 
that if the quartics in Step 2 are linearly dependent
then the invariants are necessarily zero. By Lemma~\ref{auxunique} 
we have $S(\phi) = 0$. Then $\Delta(\phi)=0$ by Lemma~\ref{auxexists}(i).
Since $c_4^3-c_6^2 = 1728 \Delta$ it only remains to show that $c_4(\phi)=0$. 
We do this by constructing a covariant
  \[ T  : \wedge ^2 V \otimes W \to S^5W^* \]
of degree 30 with $\langle S,T \rangle = c_4^2$. 
We omit the (lengthy) details, since our main interest is in applying
Algorithm~\ref{myalg} in the case $C_{\phi}$ 
is a smooth curve of genus one.

An alternative method for computing the discriminant is the following.
Let $\phi \in X_5$ with $4 \times 4$ Pfaffians 
$p_1, \ldots, p_5$. For a permutation $\pi \in S_5$ we define
\[ \begin{array}{rcl} \Omega_{\pi(1),\pi(2)} & = & \sign(\pi) 
\sum_{i,j=1}^5 \frac{\partial p_i}{\partial x_{\pi(3)}} 
\frac{ \partial \phi_{ij}}{\partial x_{\pi(4)}}   
\frac{\partial p_j}{\partial x_{\pi(5)}}.  \end{array} \]
The calculations of \S\ref{sec:invdiff} show that $\Omega_{r,s}$ is 
well-defined up to the addition of quadrics in the space
spanned by $p_1, \ldots, p_5$.
We arrange the coefficients of $p_1, \ldots, p_5$ and $\Omega_{r,s}$
for $1 \le r < s \le 5$ in a $15 \times 15$ matrix. 
Then the determinant of this matrix is an invariant of degree 60,
and hence weight 12. We claim it is $\pm 32 \Delta$.
Since the invariants of 
% We claim that the determinant of this matrix is $\pm 32 \Delta$.
% Indeed the calculations of \S\ref{sec:invdiff} show that this
% determinant is an invariant. According to Lemma~\ref{degwt}, an 
% invariant of degree 60 has weight 12. But the invariants of
weight $12$ form a 2-dimensional vector space, 
our claim is verified by computing two (suitably chosen)
numerical examples. 
This method for computing the discriminant is in practice much faster
than using Algorithm~\ref{myalg}. 

\section{Computing the geometric invariants}
\label{compgeom}

Let $C \subset \PP^{n-1}$ be a genus one normal curve of 
degree $n \ge 3$, and let $\omega$ be an invariant differential 
on $C$, both defined over a field $K$. The geometric invariants 
$c_4$ and $c_6$ of the pair $(C,\omega)$ were defined in \S\ref{geominvar}.
We are interested in computing geometric invariants for 
the following two reasons.

\medskip

\noindent
{\bf Computing the Jacobian.}
Given equations defining a genus one normal curve $C \subset \PP^{n-1}$ 
of degree $n$, we aim to compute a Weierstrass equation for its
Jacobian. The first step is to compute an invariant differential
$\omega$ on $C$. We can do this using either the method of \S\ref{geominvar}
or the method of \S\ref{sec:invdiff}. Proposition~\ref{ginvjac}
then reduces the problem of computing the Jacobian to
that of computing the geometric invariants.

\medskip

\noindent
{\bf Minimisation.}
Let $K$ be a local field with discrete valuation
$\ord : K^* \to \Z$. If $n \le 5$ then by minimisation we mean
the task of finding an integer coefficient genus one model 
equivalent to a given one, with $\ord ( \Delta)$ minimal. 
We refer to \cite{CrSt} for a treatment of this problem in 
the case $n=2$. In general the same question can be asked provided 
we have a notion of genus one model with 
the following properties:
\begin{itemize}
\item a (non-singular) genus one model defines a pair $(C,\omega)$,
\item it is possible to decide whether a genus one model 
has integer coefficients. 
\end{itemize}
We will not discuss the possible definitions of genus
one model for $n >5$, 
but merely note that if we are to keep track of our progress in
minimising, we must be able to compute geometric invariants.

\medskip

We have compiled the following list of methods for
computing geometric invariants. By Lemma~\ref{geomwts}
we are free to rescale $\omega$ at any stage (provided
we keep track of the scalars).

\subsection{The invariants method}
We assume that $C$ has degree $n \le 5$. The first step is to 
compute a genus one model $\phi \in X_n$ with $C= C_\phi$. 
For $n \le 4$ this is trivial. For $n=5$ we use the algorithm
described in \cite{g1pf}. Then the formulae and algorithms
of \S\S\ref{sec:formulae},\ref{sec:evalalg} are used to  
compute $c_4(\phi)$ and $c_6(\phi)$. By Proposition~\ref{invginv} 
these are the geometric invariants of $(C_\phi,\omega_\phi)$.
% If the invariant differential $\omega_\phi$ differs from that
% originally specified, then we use Lemma~\ref{geomwts} to 
% compensate for this in our final answer.

The main disadvantage of the invariants method 
is that we are currently restricted to $n \le 5$.

\subsection{The projection method}
Extending our field (if necessary) we first 
find a rational point $P \in C(K)$. For instance we might
find $P$ by intersecting our curve with a random hyperplane,
or by taking the generic point defined over the function field.
% For instance we might take the intersection of the curve 
% with a random hyperplane.
Then we project away from $P$ to obtain a genus one normal curve
$C_P \subset \PP^{n-2}$ of degree $n-1$. 
Explicitly, we change co-ordinates on $\PP^{n-1}$ so that 
$P=(0:0: \ldots :0:1)$ and the tangent line at $P$ 
is $x_1= \ldots =x_{n-2}=0$.  Then the projection
map \[\PP^{n-1} \ratto \PP^{n-2} \, ; \quad (x_1 : \ldots :x_{n-1}: x_n)
\mapsto (x_1 : \ldots :x_{n-1})\] restricts to an isomorphism
$\pi : C \isom C_P$ with $\pi(P)= (0:0: \ldots : 0:1)$. 
We eliminate $x_n$ from the quadrics generating $I(C)$ 
by linear algebra. If $n \ge 5$ then by Proposition~\ref{nquads}
the remaining quadrics are sufficient to generate $I(C_P)$. 
The invariant differential $\omega$ on $C$ is specified by an
$n \times n$ matrix of quadrics, as described in \S\ref{geominvar}. 
The corresponding invariant
differential on $C_P$ is obtained by deleting the last row 
and column of this matrix. We eliminate $x_n$ from the remaining entries 
by subtracting suitable elements of $I(C)$.

At this stage we may either project away from $\pi(P)$ or switch
to another method. If we keep projecting away from a rational
point, then eventually we obtain a curve in Weierstrass form.
(The final stages of this process are described 
in \cite[\S8]{CaL}.) Alternatively if a method for computing Riemann-Roch
spaces is available, then we may pass directly to a Weierstrass
equation by computing $\CL(mP)$ for $m=1,2,3$. 

The main disadvantage of the projection method is that it 
requires a field extension. 

\subsection{The covering method}
Suppose we are given pairs $(C_1,\omega_1)$ and $(C_2,\omega_2)$,
and a morphism $\pi : C_1 \to C_2$. 
Further suppose that $\pi$ is a 
twist of the multiplication-by-$m$ map on an elliptic curve. 
Then $(C_1, \pi^* \omega_2)$ and $(C_2,m \omega_2)$
have the same geometric invariants. 
This enables us to compute the geometric invariants 
of $(C_1,\omega_1)$ from those of $(C_2,\omega_2)$. 

The main disadvantage of the covering method is that we need to know a 
suitable map $\pi: C_1 \to C_2$. 
However if the curve $C_1$ is found by a descent 
calculation then it is likely that such a map will be known.
In this setting we already know the Jacobian, and the
application we have in mind is minimisation.

\subsection{The Wronskian method} 
The invariant differential $\omega$ determines a derivation $f \mapsto 
d f/ \omega$ on the function field $K(C)$. Anderson \cite{And}
gives a formula in terms of Wronskian determinants for the covering
map of degree $n^2$ from $C$ to its Jacobian. 
From this data it is easy to read off the geometric invariants.

The main disadvantage of the Wronskian method is that 
it requires extensive calculations in the function field.

\medskip

\noindent
{\bf An example.}
Wuthrich \cite{Chris} has constructed an element of order $5$ in 
the Tate-Shafarevich group of an elliptic curve $E$ over $\Q$,
where the elliptic curve $E$ does not admit any rational 5-isogenies.
Written as a genus one normal quintic his example has equations
\[ \begin{array}{lcl}
p_1 & = & 3 x_1^2+x_1 x_5-x_2 x_4-x_3^2 \\
p_2 & = & 17 x_1^2-10 x_1 x_3+7 x_1 x_5-7 x_2 x_4-4 x_2 x_5+4 x_3 x_4 \\
p_3 & = & 215 x_1^2-16 x_1 x_2-80 x_1 x_3+16 x_1 x_4+81 x_1 x_5-49 x_2 x_4  \\
& & \multicolumn{1}{l}{\quad ~-28 x_2 x_5-16 x_3 x_5-16 x_4^2} \\
p_4 & = &  60 x_1^2+48 x_1 x_2-34 x_1 x_3-24 x_1 x_4+20 x_1 x_5-8 x_2^2-5 x_2 x_3  \\
& & \multicolumn{1}{l}{\quad ~-12 x_2 x_4+16 x_2 x_5-14 x_3 x_5-8 x_4 x_5} \\
p_5 & = & 18 x_1^2+9 x_1 x_3-4 x_1 x_4-4 x_1 x_5-4 x_2 x_3-8 x_2 x_4-6 x_2 x_5 \\ & & \multicolumn{1}{l}{\quad ~+8 x_3 x_5-4 x_5^2.}
\end{array} \]
We use the algorithm in \cite{g1pf} to write these quadrics as the
$4 \times 4$ Pfaffians of a matrix of linear forms:
\tiny 
\[ \begin{pmatrix}
    0 &
    310 x_1 + 3 x_2 + 162 x_5 & 
    -34 x_1 - 5 x_2 - 14 x_5 & 
    10 x_1 + 28 x_4 + 16 x_5 & 
    80 x_1 - 32 x_4  \\ 
%   -310 x_1 - 3 x_2 - 162 x_5 
    & 0 & 
    6 x_1 + 3 x_2 + 2 x_5 & 
    -6 x_1 + 7 x_3 - 4 x_4 & 
    -14 x_2 - 8 x_3 \\
%     34 x_1 + 5 x_2 + 14 x_5 & 
%     -6 x_1 - 3 x_2 - 2 x_5 & 
    & & 0 & 
    -x_3 & 
    2 x_2 \\ 
%    -10 x_1 - 28 x_4 - 16 x_5 & 
%    6 x_1 - 7 x_3 + 4 x_4 & 
%    x_3 & 
    & - & & 0 & 
    -4 x_1 \\ 
%    -80 x_1 + 32 x_4 & 
%    14 x_2 + 8 x_3 & 
%    -2 x_2 & 
%    4 x_1 & 
    & & & & 0
\end{pmatrix} \]
\normalsize
Algorithm~\ref{myalg} then computes the invariants
 \[c_4 = 2^{44} \times 151009,  % 2656578422381215744,
\qquad c_6 = -2^{66} \times 34871057.  \] 
Thus the Jacobian is the elliptic curve of conductor
$1\,289\,106\,508\,910$ with minimal Weierstrass equation
\[ y^2+xy+y=x^3+x^2-3146 x +39049.  \]
% It has square-free conductor $N=1289106508910$.
According to {\sf MAGMA} \cite{magma} this elliptic curve has rank 0 
and the analytic order of its Tate-Shafarevich group is 25.
It is also the only elliptic curve in its isogeny class.

\medskip

We were also able to compute this example using the
projection and Wronskian methods. 
% (The covering method does not apply.) 
In our current implementation (written in {\sf MAGMA} \cite{magma},
and available from the author's website) the invariants method
is slightly faster than the projection method, each taking around 
a second. The Wronskian method is much slower, taking around 30
seconds in this case, but has the 
advantage of giving equations for the covering map.
These timings are of course heavily dependent on details of
the implementation we have not described here.

\section{Invariants in characteristics $2$ and $3$}
\label{sec:char23}

In \S\ref{invars} we showed that there is an injective homomorphism
of graded rings 
\[ \pi_n^* : K[X_n]^{G_n} \to K[X_1]^{G_1}. \]
We also recalled the usual formulae for $b_2,b_4,b_6,b_8$ and
$c_4,c_6,\Delta$ as polynomials in \[K[X_1] = K[a_1,a_2,a_3,a_4,a_6].\] 
In Lemma~\ref{inv1} we saw that if $\charic(K)  \not=2,3$ 
then $K[X_1]^{G_1} = K[c_4,c_6]$.
The analogue of this result in characteristics $2$ and $3$ is the following.
% We now give the analogous result in characteristics $2$ and $3$.
\begin{Lemma} The ring of invariants is % If $\charic(K) = 2$ or $3$ then
\[ K[X_1]^{G_1} = \left\{ \begin{array}{ll} 
K[a_1,\Delta] & \text{ if } \charic(K) = 2 \\
K[b_2,\Delta] & \text{ if } \charic(K) = 3. 
\end{array} \right. \]
\end{Lemma}
\begin{Proof} 
It is easy to show that $a_1$ and $\Delta$, 
respectively $b_2$ and $\Delta$, are invariants. 
We must show that they generate the ring of invariants.
As in the proof of Lemma~\ref{inv1}, this is deduced from
the existence of a suitable normal form.

\paragraph{Case $\charic(K)= 2$.} We start with the general Weierstrass
equation 
\[ y^2 + a_1 xy + a_3 y = x^3 + a_2 x^2 + a_4 x + a_6. \]
One easily computes $j=a_1^{12}/\Delta$. We assume $j \not=0$ and
following \cite[Appendix A]{Si1} make substitutions 
$x= x' +r$ and $y=y'+t$ so that $a_3=a_4=0$. We are free to suppose that
$K$ is algebraically closed. Then a further substitution $y=y'+sx'$ 
gives $a_2=0$. We arrive at the normal form
\[ y^2 + a_1 xy  = x^3 + a_6 \]
with $\Delta = a_1^6 a_6$. It follows that every invariant
is a polynomial in $a_1$ and $\Delta/a_1^6$. We are done since
$a_1$ does not divide $\Delta$. 

\paragraph{Case $\charic(K)= 3$.} We start with a general Weierstrass
equation and complete the square to obtain
\[ y^2 = x^3 + a_2 x^2 + a_4 x + a_6. \]
One easily computes $j=a_2^6/\Delta$. We assume $j \not=0$ and
following \cite[Appendix A]{Si1} make a substitution 
$x= x' +r$ so that $a_4=0$. We arrive at the normal form
\[ y^2 = x^3 + a_2 x^2 + a_6 \]
with $b_2 = a_2$ and $\Delta = -a_2^3 a_6$. It follows that every 
invariant is a polynomial in $b_2$ and $\Delta/b_2^3$. We are done since
$b_2$ does not divide $\Delta$. 
\end{Proof}

\begin{Theorem} Let $n=2,3,4,5$. Then the 
map $\pi_n^* : K[X_n]^{G_n} \to K[X_1]^{G_1} $
is an isomorphism in all characteristics.
\end{Theorem}
\begin{Proof}
In \S\ref{invars} we saw that $c_4,c_6,\Delta \in  K[X_1]^{G_1}$
extend to invariants $c_4,c_6,\Delta \in  K[X_n]^{G_n}$.
So it only remains to show that in characteristics~$2$ and~$3$
there are invariants in $K[X_n]^{G_n}$ of weights 1 and 2.

If $\charic(K) = 2$ or $3$ then $c_4^3 - c_6^2 = 1728 \Delta = 0$.
So an invariant of weight 2 exists by unique factorization in $K[X_n]$.

We now take $\charic(K)=2$ and split into the cases $n=2,3,4,5$.
% It only remains to show that if $\charic(K)=2$ then $K[X_n]^{G_n}$
% contains an invariant of weight 1. 
In the cases $n=2,3$ the coefficient of $xyz$
is an invariant of weight 1. In the case $n=4$ we write
\[ \begin{array}{rcl}
 q_1 (x_1, \ldots ,x_4) & = & \sum_{i \le j} a_{ij} x_i x_j \\
 q_2 (x_1, \ldots ,x_4) & = & \sum_{i \le j} b_{ij} x_i x_j 
\end{array} \]
and find % that the required invariant is
\[ a_1 = a_{12} b_{34} + a_{13} b_{24} + a_{14} b_{23} 
 + a_{23} b_{14} + a_{24} b_{13} + a_{34} b_{12}. \]
If $n=5$ then our genus one model is a matrix of linear
forms, say $\phi = (\phi_{ij}(x_1, \ldots,x_5))$. Let $T$ be a set of left 
coset representatives for $D_{5} = \langle (12345),(25)(34) \rangle$ 
as a subgroup of $S_5$. 
Then $a_1$ is the coefficient of $\prod_{i=1}^5 x_i$ in  
$ \sum_{\sigma \in T} \prod_{i=1}^5 \phi_{\sigma(i) \, \sigma(i+1)}$.
\end{Proof}

\begin{Remark}
If $\charic(K) =2$ or $3$ then the invariants do not suffice to 
compute the Jacobian. For example the elliptic curves
$y^2+xy= x^3+1$ and $y^2+xy= x^3+x^2+1$ over $\F_2$ 
have invariants $a_1= \Delta=1$, but are not isomorphic. 
Similarly the elliptic curves $y^2=x^3-x \pm 1$ over $\F_3$
have invariants $b_2=0$ and $\Delta=1$, but are not isomorphic.
These examples should be seen as a consequence of the failure
of Lemma~\ref{om1} in characteristics $2$ and $3$.
% the case $\charic(K) =2$ or $3$. 
\end{Remark}

As we noted in the introduction, it should instead 
be possible to find a formula for the Jacobian that works 
in all characteristics by modifying the formulae in characteristic $0$.
This has been carried out by Artin, Rodriguez-Villegas and Tate 
\cite{ARVT} in the case $n=3$.

\subsection*{Acknowledgements}
I would like to thank Nick Shepherd-Barron for
introducing me to this problem.
The computer calculations in support of this work were
performed using {\sf MAGMA} \cite{magma} and {\sf PARI} \cite{pari}.


\begin{thebibliography}{99}

\bibitem%[AKMMMP]
{Mc+}
S.Y. An, S.Y. Kim, D.C. Marshall, S.H. Marshall,
W.G. McCallum and A.R. Perlis,
Jacobians of genus one curves,
{\em J. Number Theory} {\bf{90}} (2001), no. 2, 304--315. 

\bibitem%[And]
{And}
G.W. Anderson, 
Lacunary Wronskians on genus one curves, 
{\em J. Number Theory} {\bf{115}} (2005), no. 2, 197--214.

\bibitem{ARVT}
M. Artin, F. Rodriguez-Villegas and J. Tate, 
On the Jacobians of plane cubics,
{\em Adv. Math.} {\bf{198}} (2005), no. 1, 366--382.

\bibitem%[BH]
{CMrings}
W. Bruns and J. Herzog, 
{\em Cohen-Macaulay rings},
Cambridge Studies in Advanced Mathematics {\bf{39}},
Cambridge University Press, Cambridge, 1993. 

\bibitem%[BE1]
{BE1}
D.A. Buchsbaum and D. Eisenbud, 
Gorenstein ideals of height $3$,
{\em Seminar D. Eisenbud/B. Singh/W. Vogel}, Vol. 2, 
pp. 30--48, Teubner-Texte zur Math., {\bf{48}}, Teubner, Leipzig, 1982. 

\bibitem%[BE2]
{BE2}
D.A. Buchsbaum and D. Eisenbud,
Algebra structures for finite free resolutions, and some 
structure theorems for ideals of codimension 3,
{\em Amer. J. Math.} 
{\bf{99}} (1977) 447-485.

\bibitem
{CaL}
J.W.S. Cassels, 
{\em Lectures on elliptic curves}, CUP, Cambridge, 1991. 

\bibitem{CrSt}
J.E. Cremona, M. Stoll, 
Minimal models for 2-coverings of elliptic curves,
{\em LMS J. Comput. Math.} {\bf{5}} (2002), 220--243.

\bibitem%[D]
{Dolg}
I. Dolgachev, {\em Lectures on invariant theory},
LMS Lecture Note Series  {\bf{296}},
Cambridge University Press, Cambridge, 2003. 

\bibitem%[E]
{E}
D. Eisenbud, 
{\em Commutative algebra with a view toward algebraic geometry},
GTM {\bf{150}}, Springer-Verlag, New York, 1995. 

\bibitem%
{hsecenc} 
T.A. Fisher,
{\em The higher secant varieties of an elliptic normal curve},
preprint.

\bibitem%
{g1pf} 
T.A. Fisher,
{\em Genus one curves defined by Pfaffians}, 
preprint.

\bibitem%[FH]
{FH} 
W. Fulton and J. Harris,
{\em Representation theory},
GTM {\bf{129}},
Springer-Verlag, New York, 1991. 

\bibitem%[GY]
{GY} 
J.H. Grace and A. Young,
{\em The algebra of invariants}, 
Cambridge University Press, Cambridge, 1903. 

\bibitem%[Ha]
{Ha}
R. Hartshorne, 
{\em Algebraic geometry}, GTM {\bf{52}}, 
Springer-Verlag, New York-Heidelberg, 1977.

\bibitem{Hilbert}
D. Hilbert, 
{\em Theory of algebraic invariants},
Cambridge University Press, Cambridge, 1993. 

\bibitem%[Hu]
{Hu} 
K. Hulek,
{\em Projective geometry of elliptic curves},
Soc. Math. de France, Ast\'erisque
{\bf{137}} (1986).

\bibitem{Kempf}
G. Kempf,
Some quotient surfaces are smooth,
{\em Michigan Math. J.} {\bf{27}} (1980), no. 3, 295--299.

\bibitem{magma} 
{\sf MAGMA} is described in W. Bosma, J. Cannon and C. Playoust, 
The Magma algebra system I: The user language, {\em J. Symb. Comb.} {\bf{24}}, 
235-265 (1997). (See also the Magma home page at 
{\tt http://magma.maths.usyd.edu.au/magma/}.)

\bibitem%[MSS]
{MSS}
J.R. Merriman, S. Siksek, N.P. Smart,
Explicit $4$-descents on an elliptic curve, 
{\em Acta Arith.} {\bf{77}} (1996), no. 4, 385--404. 

\bibitem{pari}
PARI/GP is developped by the PARI~Group, University of Bordeaux.
(See also the PARI home page at 
{\tt http://pari.math.u-bordeaux.fr/}.)

\bibitem%[Sa]
{Salmon} 
G. Salmon, 
{\em A treatise on the higher plane curves},
Third edition, Hodges, Foster and Figgis, Dublin, 1879.

\bibitem%[Sh]
{Sh}
I.R. Shafarevich, 
{\em Basic algebraic geometry. 1. Varieties in projective space}, 
Springer-Verlag, Berlin, 1994. 

\bibitem%[Si1]
{Si1}
J.H. Silverman, 
{\em The arithmetic of elliptic curves},
GTM {\bf{106}}, Springer-Verlag, New York, 1986. 

\bibitem%[St]
{St}
B. Sturmfels, 
{\em Algorithms in invariant theory},
Texts and Monographs in Symbolic Computation,
Springer-Verlag, Vienna, 1993. 

\bibitem%[We]
{Weil1954}
A. Weil, 
Remarques sur un m\'emoire d'Hermite,
{\em Arch. Math.} {\bf{5}} (1954), 197--202. 

\bibitem%[Wu]
{Chris}
C. Wuthrich, 
Une quintique de genre 1 qui contredit le principe de Hasse,
{\em Enseign. Math.} (2) 47 (2001), no. 1-2, 161--172.

\end{thebibliography}
\end{document}